\documentclass{article}%
\usepackage{amsfonts}
\usepackage{graphicx}
\usepackage{amsmath}
\usepackage{amssymb}%
\usepackage{float}
\providecommand{\U}[1]{\protect\rule{.1in}{.1in}}
\newtheorem {theorem}{Theorem}[section]

\newtheorem {corollary}{Corollary}[section]

\newtheorem{definition}{Definition}[section]

\newtheorem{remark}{Remark}[section]

\newtheorem {condition}{Condition}[section]

\begin{document}

\begin{center}
{\LARGE Memory properties of  transformations of linear processes}

\bigskip

\centerline{\today}
\bigskip  Hailin Sang and Yongli Sang

\bigskip Department of Mathematics, University of Mississippi,
University, MS 38677, USA. E-mail address: sang@olemiss.edu, ysang@go.olemiss.edu

\end{center}

\bigskip\textbf{Abbreviated Title: }{\Large Transformations of linear processes}

\begin{center}
\bigskip\textbf{Abstract}

\end{center}
In this paper, we study the memory properties of  transformations of linear processes. Dittmann and Granger (2002) studied the polynomial  transformations of Gaussian FARIMA$(0,d,0)$ processes by applying the orthonormality of the Hermite polynomials under the measure for the standard normal distribution. Nevertheless, the orthogonality does not hold for transformations of non-Gaussian linear processes. Instead, we use  the decomposition developed by Ho and Hsing (1996, 1997) to study the memory properties of nonlinear transformations of linear processes,  which include the FARIMA$(p,d,q)$ processes, and obtain consistent results as in the Gaussian case. In particular,  for stationary processes, the transformations of short-memory time series still have short-memory and the transformation of long-memory time series may have different weaker memory parameters which depend on the power rank of the transformation. On the other hand, the memory properties of transformations of non-stationary time series may not depend on the power ranks of the transformations. This study has application in econometrics and financial data analysis when the time series observations have non-Gaussian heavy tails.  As an example, the memory properties of call option processes at different strike prices are discussed in details. \\          

\noindent Key words and phrases:  heavy tail, long memory, linear process, nonlinear transformation, non-stationary, short memory\\

\noindent  {\textit{MSC 2010 subject classification}: 62M10, 62E20}


\section{Introduction}
Let 
\begin{equation}\label{lp}
X_n=\sum_{i=0}^{\infty}{a_i\varepsilon_{n-i}}, \;\;n\in\mathbb{N},
\end{equation}
 be a linear process, where the innovations $\varepsilon_i$, $i\in \mathbb{Z}$,  are 
 i.i.d. random variables with mean zero and  finite variances. Without loss of generality, we assume that $\mathbb{E}\varepsilon_i^2=1$, $i\in \mathbb{Z}$. The coefficients $a_i$ satisfy $\sum_{i=0}^{\infty}a_i^2<\infty$, under which the linear process (\ref{lp}) is well defined by the three series theorem.  
 
A stationary time series $X_n$ has short or long memory (short or long range dependence) in the covariance sense depending on $\sum_{n=1}^\infty |\mbox{Cov}(X_1, X_n)|<\infty$ or $=\infty$ (Parzen, 1981).
Meanwhile,  in the frequency domain, a stationary time series $X_n$ with a spectral density function $f(\lambda)$ is called a long memory process in a restricted spectral density sense if $f(\lambda)$ is bounded on $[\delta, \pi]$ for every $\delta>0$,  and  $f(\lambda)\rightarrow\infty$ as $\lambda\rightarrow 0^+$. These two definitions are not always equivalent,  
 see Cox (1977) and Gu\a'{e}gan (2005). In particular, for $d<1/2$, Dittmann and Granger (2002) called a stationary time series $X_n\sim$ LM$(d)$ if the spectral density function $f(\lambda)$ behaves like a power function at low frequencies, that is as $|\lambda|^{-2d}$ as $\lambda$ approaches zero. For $d\ge 1/2$, $X_n\sim$ LM$(d)$ if and only if $(1-B)^kX_n\sim$ LM$(d-k)$ for $k=[d+1/2]$, where $[x]$ denotes
the largest integer smaller or equal to $x$. $B$ is the backward shift operator, $BX_{i}=X_{i-1}$. The cases $d>0$, $d=0$ and $d<0$ correspond to long memory, short memory and negative dependence (antipersistence),  respectively. 
 
 Dittmann and Granger (2002) studied the memory properties of polynomial  transformation of Gaussian FARIMA$(0,d,0)$ processes
\begin{equation}
X_{n}=(1-B)^{-d}\varepsilon_{n}=\sum_{i\geq0}a_{i}\varepsilon_{n-i},\label{deffractlin}
\end{equation}
where $a_{i}=\frac{\Gamma(i+d)}{\Gamma(d)\Gamma(i+1)}$, $\varepsilon_i\sim i.i.d. \;\;N(0, \sigma^2)$, $-1<d<1/2$ and $d\ne 0$. $d=0$ gives the i.i.d. process $\{\varepsilon_n\}$.  They applied the orthonormality of the Hermite polynomials under the measure for the standard normal distribution. That is, 
\begin{align}\label{ort}
\int_{-\infty}^\infty H_m(x)H_n(x)dP(Z\le x)=I(m=n),
\end{align}
where Hermite polynomials $H_j(x)$ are defined by 
$$\left(\frac{d}{dx}\right)^j e^{-x^2/2}=(-1)^j\sqrt{j!}H_j(x)e^{-x^2/2},  j=0,1,2, \cdots, $$  $Z\sim N(0,1)$ and $I(m=n)$ is the indicator function, $m, n=0,1,2, \cdots$. For example, see Cram\'{e}r (1946). In the continuous case, Taqqu (1979) and Giraitis and Surgailis (1985) studied the nonlinear transformations of fractional Brownian motions. Nevertheless, this nice orthogonal property (\ref{ort}) does not hold in general when the distribution is not Gaussian. On the other hand, it is witnessed and well known in the financial field that quite many financial data like stock prices have heavier tails than the tail of the normal distribution, for example, see Ruppert (2011). For the non-Gaussian case, based on the innovations $\varepsilon_i$, Ho and Hsing (1996, 1997) developed an expansion with orthogonal terms which is akin to the Hermite expansion for the Gaussian case. 

We focus on the transformations of linear processes (\ref{lp}), which are not necessarily Gaussian, in this paper. We assume that $a_0=1$, $\sum_{i=1}^\infty|a_i|<\infty$ or 
$a_i=i^{-\beta}L(i)$, $i>0$,  for some $\beta \in (1/2, 1)$,
 where $L(i)>0$ is a slowly varying function at $\infty$ (Bingham, Goldie and Teugels, 1987), i.e., $\lim_{x\rightarrow\infty}L(\lambda x)/L(x)=1$ for any $\lambda>0$. It includes the autoregressive fractionally integrated moving average FARIMA$(p,d,q)$ processes introduced by Granger and Joyeux (1980) and Hosking (1981), which is defined as
\begin{align}\label{farima}
\phi(B)X_n=\theta(B)(1-B)^{-d}\varepsilon_{n}.
\end{align}
 Here $p, q$ are nonnegative integers, $\phi(z)=1-\phi_1z-\cdots -\phi_pz^p$ is the AR polynomial and $\theta(z)=1+\theta_1z+\cdots \theta_qz^q$ is the MA polynomial. Under the conditions that $\phi(z)$ and $\theta(z)$ have no common zeros, the zeros of $\phi(\cdot)$ lie outside the closed unit disk and $-1<d<1/2$, the FARIMA($p,d,q$) process has linear process form (\ref{lp}) with $a_i=\frac{\theta(1)}{\phi(1)}\frac{i^{d-1}}{\Gamma(d)}+\mathcal{O}(i^{-1})$. See Bondon and Palma (2007) for the extension of causality to the range of $-1<d<1/2$  and Kokoszka and Taqqu (1995) for the asymptotic coefficient formula. 
 
To explore the  memory properties of transformations $K(X_n)$ with $\mathbb{E}K^2(X_n)<\infty$ of linear processes (\ref{lp}),  we shall apply the decomposition of $K(X_n)$ proposed by Ho and Hsing (1996, 1997). See also the review paper by Hsing (2000). This method has been applied to the expansion of $K(X_n)$ for linear processes $X_n$ in the study of many subjects, for examples, weak convergence theorems including central limit theorem, functional central limit theorem, convergence to Wiener-Ito integral and Hermite process (Ho and Hsing, 1996, 1997; Hsing, 1999; Wu, 2002, 2006), the kernel density estimation (Honda, 2000; Wu and Mielniczuk, 2002; Kulik, 2008), the empirical processes of long memory sequences (Wu 2003),  the U-statistics (Ho and Hsing, 2003; Hsing and Wu, 2004) and the moderate deviations (Wu and Zhao, 2008; Peligrad et al., 2014).  Under the condition proposed by Wu (2006), we obtain results both in time domain and frequency domain. The results in time domain is consistent with the limit theorems in Ho and Hsing (1997) and Wu (2006) via the order of normalization. The results are applicable not only to FARIMA$(0,d,0)$ processes as studied in Dittmann and Granger (2002) for Gaussian case, but also to general FARIMA$(p,d,q)$ processes for some special transformation. The results hold not only for smooth transformations, they also hold for functions which are not differentiable. In particular, we study the memory properties of option time series $(X_n-C)^+$ in finance for different strike price $C>0$. 

We also study the  properties of nonlinear transformations of non-stationary time series $X_n$ with the form 
\begin{align}\label{nons1}
X_n=\sum_{j=1}^n Y_j,  \;\; \text{where} \;\;\;Y_j=\sum_{i=0}^{\infty}{a_i\varepsilon_{j-i}},
\end{align}
$a_{i}=\frac{\Gamma(i+d-1)}{\Gamma(d-1)\Gamma(i+1)}$, $1/2<d<1$.
Again, we do not assume that the innovations are Gaussian.       

In this paper we shall use the following notations:  we use  $a_{m}\mathbb{\sim
}b_{m}$ instead of the notation $a_{m}/b_{m}\rightarrow1$; for positive sequences, the
notation $a_{m}\ll b_{m}$ or $b_{m}\gg a_{m}$ and the Vinogradov symbol $\mathcal{O}$ mean that
$a_{m}/b_{m}$ is bounded; the notation $a_m\simeq b_m$ means that there exist constants $c_1$ an $c_2$ such that $0<c_1b_m<a_m<c_2b_m$ for $m$ large enough.  $C>0$ is a generic constant which may vary in different context. 

The paper has the following structure. In Section \ref{main} we study the memory properties of transformations of the stationary long-memory and short-memory processes. Section 3 is on the non-stationary processes. In Section 4, we study the application to option processes in finance. The proofs go to Section 5. 


\section{ Transformations of stationary  processes}\label{main}
In this section we study the  transformation $K(X_n)$ of the stationary process (\ref{lp}). 
First of all, $K(X_n)$, $n\in\mathbb{N}$, is strictly stationary since the time series $X_n$ is strictly stationary. By the condition $\mathbb{E}K^2(X_n)<\infty$, $K(X_n)$ is also (covariance) stationary. We start this section with a couple of basic notations which will be used throughout the paper.
Let $||X||=[\mathbb{E}(X^2)]^{1/2}$ be the $L^2$ norm of the random variable $X.$
Define the shift process  $\mathcal{F}_i=(\ldotp \ldotp \ldotp,\varepsilon_{i-1}, \varepsilon_i  )$,  and let
%
\begin{align*}
& X_{n,k}=\mathbb{E}(X_n|\mathcal{F}_k),\\
& K_n(w)=\mathbb{E}[K(w+X_n-X_{n,0})],\\
&K_\infty(w)=\mathbb{E}[K(w+X_n)],\\
&K^{(r)}_n(w)=\frac{d^r}{d w^r}\mathbb{E}[K(w+X_n-X_{n,0})] ,\\
&K^{(r)}_\infty(w)=\frac{d^r}{d w^r}\mathbb{E}[K(w+X_n)]
\end{align*}    
for any nonnegative integer $r$. 
The following definition is from Ho and Hsing (1997). 
\begin{definition}\label{def}
 A transformation $K(\cdot)$ has power rank $k$ with respect to the linear process $X_n$ for some positive integer $k$ if $K_{\infty}^{(k)}(0)\ne 0$  and $K_{\infty}^{(r)}(0)=0$ for all $1 \leq r <  k$.
\end{definition}
We will use $k$ to denote the power rank of $K(\cdot)$ with respect to the linear process $X_n$ throughout the paper. 

In order to present the main results of this paper, we need the following condition from Wu (2006). 
\begin{condition}\label{cond}
Let $\mathbb{E}(|\varepsilon_1|^q)<\infty$ for some $2<q\leq 4$ and $K_{n} \in \mathbb{C}^{k+1}(\mathbb{R})$ for all large $n.$ Assume that for some $\lambda>0,$
\begin{align*}
\sum_{\alpha=0}^{k+1}\|K^{(\alpha)}_{n-1}(X_{n,0};\lambda)\|+\sum_{\alpha=0}^{k-1}\||\varepsilon_1|^{q/2}K^{(\alpha)}_{n-1}(X_{n,1})\|+\|\varepsilon_1K^{(k)}_{n-1}(X_{n,1})\|=\mathcal{O}(1),
\end{align*}
\end{condition}
where  $K^{(\alpha)}_{n-1}(X_{n,0};\lambda)=\mbox{sup}_{|y| \leq \lambda}|K^{(\alpha)}_{n-1}(X_{n,0}+y)|$ is the  local maximal function for $K^{(\alpha)}_{n-1}(X_{n,0})$.  As Wu (2006) mentioned,  Condition \ref{cond} is quite mild,  which only imposes certain smoothness requirements on $K_{n-1}.$

First we consider the case $a_0=1$ and $ a_i=i^{-\beta}L(i), i>0$, $1/2<\beta<1$,  for the linear process (\ref{lp}). Notice that in this case the covariance function $\gamma_X(h)=\mathbb{E}X_0X_h$ of the original series $X_n$ is regularly varying with exponent $-1<1-2\beta<0$ and hence $X_n$ has long memory in the covariance sense. The FARIMA$(p,d,q)$ process as in (\ref{farima}) with $0<d<1/2$ is a particular example of this case, $\beta=1-d$. 
\begin{theorem}\label{thm1}
Assume that Condition \ref{cond} holds with $q=4$ and that $K$ has power rank $k\geq 1.$ Let $a_0=1$ and $ a_i=i^{-\beta}L(i), i>0$, $1/2<\beta<1$, in model (\ref{lp}).  If the power rank $k$ of a transformation $K(\cdot)$ with respect to the linear process (\ref{lp}) satisfies $k<(2\beta-1)^{-1}$, then $K(X_n)$ has long memory in the covariance sense. $K(X_n)$ has short memory in the covariance sense if $k>(2\beta-1)^{-1}$. 
\end{theorem}
This theorem shows that $K(X_n)$ has long memory as long as the power rank of $K(\cdot)$ satisfies $k<(2\beta-1)^{-1}$. Hence $K(X_n)$ keeps the long memory property for a wide range (in terms of the power rank $k$) of transformations if the parameter $\beta$ of the original series $X_n$ is close to $1/2$ and therefore $X_n$ has very strong long memory. 
Nevertheless, $K(X_n)$ losses the long memory property for a wide range of transformations if  $\beta$ is not close to $1/2$. For example, if $3/4<\beta<1$, only $X_n$ and other transformations with power rank $k=1$ keep the long memory property. 
\begin{remark}
Ho and Hsing (1996, 1997) studied the limit theorems by assuming  $\mathbb{E}K^2(X_n)<\infty$, $\mathbb{E}\varepsilon^8<\infty$ and the condition $C(t, \tau, \lambda)$ there. Wu (2006) studied the functional limit theorems under the improved condition $\mathbb{E}\varepsilon^4<\infty$ and the Condition \ref{cond}. The memory property in Theorem  \ref{thm1} is consistent with the above limit theorems via the order of normalization.  Theorem  \ref{thm1} only requires the Condition \ref{cond} with $q=4$.  
\end{remark}
\begin{remark}
It is well known that both long memory and  heavy tail parameters play roles in the asymptotics of partial sums of time series. See e.g., McElroy and Politis (2007). As in Wu (2006), for the major results in this paper, we assume that the innovation of the linear process has fourth moment, although not necessary Gaussian. Therefore, the memory parameter dominates the growth of the partial sum, which is the case in the upper right hand region of Figure 1 in McElroy and Politis (2007).  
\end{remark}
The next corollary shows that if the slowly varying function $L(x)$ is a constant asymptotically, $K(X_n)$ has long memory also in the case that $(2\beta-1)^{-1}$ is an integer and the power rank $k=(2\beta-1)^{-1}$. 
\begin{corollary}\label{corollary1}
In the case that $\lim_{n\rightarrow \infty}L(n)=L$ for some constant $L>0$, under the same conditions as in Theorem \ref{thm1},  $K(X_n)$ has long memory in the covariance sense if $k\le(2\beta-1)^{-1}$. $K(X_n)$ has short memory in the covariance sense if $k>(2\beta-1)^{-1}$. 
\end{corollary}
Corollary \ref{corollary1} is applicable to FARIMA$(p,d,q)$ process. Recall that in this case, $\beta=1-d$ and $\lim_{n\rightarrow \infty}L(n)=\frac{\theta(1)}{\phi(1)\Gamma(d)}$.  Furthermore, we have detailed knowledge on the memory parameter of $K(X_n)$ from the following Theorem \ref{cor2} if the linear process  is a FARIMA$(p,d,q)$ process. 
\begin{theorem} \label{cor2}
 Let $X_n$ be a stationary  
 FARIMA$(p,d,q)$  process (\ref{farima}) with  $0<d<1/2$  and  Condition \ref{cond} holds with $q=4.$ $K(\cdot)$ has power rank $k$  with respect to the FARIMA(p,d,q) process. Then $K(X_n)$ is a long-memory process LM$(\tilde{d})$ with $\tilde{d}= (d-1/2)k+1/2$ when $k(1-2d)<1,$ and a short-memory process LM$(0)$ if $k(1-2d)>1$ but $(k-1)(1-2d)<1.$
\end{theorem}
This theorem shows that $K(X_n)$ can never have stronger long range dependence than the original process since $\tilde{d}\le d$ and $\tilde{d}=d$ if and only if $k=1$. 


Now we study the transformations of short memory linear processes in the form of (\ref{lp}). The following theorem provides a result in the general setting $\sum_{i=0}^{\infty}|a_i|<\infty$. 
\begin{theorem}\label{thm2}
Assume  $\sum_{i=0}^{\infty}|a_i|<\infty$  in the model (\ref{lp}) and
 \begin{align}\label{cond2}
\|K_{n-1}(X_{n,1})-K_{n-1}(X_{n,0})\|=\mathcal{O}(|a_{n-1}|).
\end{align} 
Then $K(X_n)$ has short memory in the covariance sense for any transformation $K(\cdot)$ with $\mathbb{E}K^2(X_n)<\infty$. 
\end{theorem}
As Condition \ref{cond}, the condition (\ref{cond2}) is proposed by Wu (2006) and only requires certain smoothness  requirements on $K_{n-1}.$ 
Theorem \ref{thm2} shows that we can never get long memory process from transformations if the original process has short memory. 

It is an open question that whether $K(X_n)$ is a short memory LM($0$) process if $X_n$ is a stationary FARIMA$(p,d,q)$ process with $-1<d<0$.  The following Theorem \ref{thm3} gives a confirmative answer to this question for the special case, the FARIMA$(0,d,0)$ process with $-1<d<0$, if $K(x)=x^2.$ The FARIMA$(0,d,0)$ is not necessarily Gaussian. 
\begin{theorem}\label{thm3}
Let $X_n$ be a stationary  FARIMA$(0,d,0)$ process, $-1<d<0$, defined as in (\ref{deffractlin}). Then $X^{2}_n$ is a short-memory process LM($0$).
\end{theorem}
Similar to the Gaussian case, Theorem \ref{thm3} shows that antipersistence is a much more fragile property than long memory property. The antipersistence is immediately lost for the square transformation.

To verify the main results in this section, in particular Theorem \ref{cor2} and Theorem \ref{thm3},  we conduct simulation study for the memory of some common transformations of FARIMA$(p,d,q)$ processes. These transformations include $K(x)=x^2$, $x^3$, $x^4$, $x^3-3x$, $x^4-6x^2$, $\sin x$,  $e^x$ and the non-continuous indicator function $I(x\leq c)$ for some constant $c$. First of all, we calculate the power rank of $K(\cdot)$ with respect to $X_n$ and then find the theoretical memory parameter of each transformed process from Theorem \ref{cor2}. Although the power rank of $K(X)$ is identical to its Hermite rank if $X$ has standard normal distribution (Ho and Hsing, 1997), it may be different under different distributions. Nevertheless, one can easily find the power rank of a specific transformation under different distributions by the Definition \ref{def}. For example, provided that $\int \cos ydF(y)\ne 0$ or $\int e^ydF(y)<\infty$, the power rank of $K(x)=\sin x$ or $K(x)=e^x$ is 1 since 
$$K_\infty(x)=\sin x\int \cos ydF(y)+\cos x\int\sin ydF(y)$$ or 
$K_\infty(x)=e^x\int e^ydF(y)$ satisfies $K_\infty'(0)\ne 0$.    
By similar analysis as above,  the transformations $K(x)=x^2$, $x^3$, $x^4$, $x^3-3x$, $x^4-6x^2$, $\sin x$ and $e^x$ have power rank 2, 1, 2, 3, 4, 1 and 1 respectively under some regular conditions on $X_n$ (the conditions for different transformations may be different). For the indicator function $K(x)=I(x\leq c)$, the power rank depends on the value of the constant $c$ by the following argument: Let $F(x)$ be the distribution function of X and assume that the density function $f(x)$ of X exists. Then 
$$K_\infty(x)=\int I(x+y\leq c)dF(y)=F(c-x)$$  
and $K^{'}_\infty(x)=-f(c-x)$. We then have $K^{'}_\infty(0)=-f(c)\neq 0$ if $f(c)$ is quite far away from 0. Under this condition,  the power rank of the indicator function is 1. If $f(c)$ is very close to 0 but $f^{'}(c)$ exists and is quite far away from 0, we can say that the power rank of this indicator function $K(\cdot)$ is 2. Under certain smoothness condition, we can continue this procedure to find the power rank of $K(\cdot)$ for $c$ in different ranges.

Secondly, to compare with the theoretical memory parameters, we perform simulation study for these transformations of  FARIMA$(p,d,q)$ processes $X_n$ with memory parameters $d=0.2$ and $0.4$. The three processes in the simulation study are FARIMA$(0,d,0)$ process (when $K(x)=x^2$, we also consider the cases that $d=-0.8$, $-0.4$, $-0.2$), FARIMA$(1,d,0)$ process with the AR coefficient $\phi_1=-0.3$,   and FARIMA$(1,d,1)$ process with the AR coefficient $\phi_1=-0.4$ and the MA coefficient $\theta_1=0.7$. 

Since our results require $\mathbb{E}K^2(X)<\infty$, $\mathbb{E}\varepsilon^4<\infty$ and some transformations involve $x^4$, we take the Student $t$ distribution with degree freedom $10$ as the innovations of the FARIMA$(p,d,q)$ processes for all transformations in our study except the last one $K(x)=e^x$.  We choose the Gaussian FARIMA$(p,d,q)$ processes in the transformation $K(X)=e^{X}$ since $\mathbb{E}e^{2X}<\infty$ is required.  
For each of these three processes and for each $d$, we conduct $N=2,000$ simulations with $n=2,000$ observations in each process by  applying the algorithm in 
Fa\"{y} et al. (2009). The memory parameters of each process and their transformations are estimated by the Fourier regression method proposed in Geweke and Porter-Hudak (1983). As studied in Hurvich et al. (1998), we choose the bandwidth $[n^{4/5}]$ for each estimation. 
The theoretical memory parameters and the estimated values are listed in Tables \ref{tab:1} and  \ref{tab:2}  respectively for each of these three processes.  We also report the empirical standard error of the $N=2,000$ estimates for each process in these tables. When $d$ is negative, $d=-0.8, -0.4$ and $-0.2$,  the theoretical memory parameters of all transformations of FARIMA$(0,d,0)$ except the square of the FARIMA$(0,d,0)$,   are left in blank since we do not have theoretical results for  these cases. When $d=0.2,$ we need $k\le 2$ by the condition on $d$ in Theorem 2.2, so  the theoretical memory parameters of transformations with rank greater than  $2$ are left blank.

\begin{table}[H]
\center
\caption{Average estimated memory parameters of some transformations of $2,000$ simulated stationary FARIMA$(0,d,0)$ processes with  $2,000$ observations in each process and $t(10)$ innovations (except the transformation $e^x$, for which we use Gaussian innovations since $\mathbb{E}e^{2X}<\infty$ is required).}
\label{tab:1}
\bigskip
\scriptsize
\begin{tabular}{ll  |  lllll}
\hline  \hline                                                     
  $K(X)$  and its   &          &\multicolumn{ 4}{c}{Memory parameter of the original series $X$}        \\
 power rank  &    & $d=-0.8$ &  $d=-0.4$ &              $d=-0.2$ &            $d=0.2$ &     $d=0.4$         \\
\hline

 $X$             &           Theory &              -0.8&              -0.4&                       -0.2   &                    0.2    &     0.4                    \\

         (rank $1$)  &    Simulation &   $-0.7674$&  $-0.4008$ &    $-0.2005$  &  $0.2042$  &   $0.4075 $   \\
         & Std error & $0.0496$&  $0.0335$  & $0.0333$  & $0.0319$  & $0.0332$  \\

\hline

         $X^2$             &           Theory &       0 &                    0 &                       0    &  0                   &     0.3                     \\

         (rank $2$)  &    Simulation &   $0.0387$&  $0.0250$ &    $0.0094$  &  $0.0405$  &   $0.2755 $   \\
         & Std error & $0.0330$& $0.0323$  & $0.0322$  & $0.0364$  & $0.0648$  \\

\hline

         $X^3$             &           Theory &             &              &                          &                    0.2     &     0.4                     \\

         (rank $1$)  &    Simulation &   $-0.1501$&  $-0.0895$ &    $-0.0540$  &  $0.0960$  &   $0.2824$   \\
 & Std error & $0.0462$& $0.0383$  & $0.0353$  & $0.0372$  & $0.0561$  \\

\hline
 $X^4$             &           Theory &                  &         &                         &       0                  &     0.3                     \\

         (rank $2$)  &    Simulation &  $0.0330$&   $0.0144$ &    $0.0038$  &  $0.0157$  &   $0.1855$   \\
          &Std error &$0.0329$&  $0.0329$  & $0.0292$  & $0.0360$  & $0.0790$  \\
          
\hline

         $X^3-3X$             &           Theory &             &              &                           &                     &     0.2                     \\

         (rank $3$)  &    Simulation &  $-0.0757$&    $-0.0160$ &    $-0.0029$  &  $0.0087 $  &   $0.2049 $   \\
 &Std error &$0.0481$&  $0.0345$  & $0.0321$  & $0.0347$  & $0.0800$  \\
\hline

         $X^4-6X^2$             &           Theory &                  &           &                          &                         &     0.1                    \\

         (rank $4$)  &    Simulation &  $0.0257$&    $0.0051$ &    $0.0020 $  &  $0.0008 $  &   $0.1138$   \\
          &Std error &$0.0301$&  $0.0294$  & $0.0317$  & $0.0322$  & $0.0882$  \\

\hline

         $\sin X$             &           Theory &                &           &                          &                    0.2     &     0.4                     \\

         (rank $1$)  &    Simulation &   $-0.1651$&   $-0.1863$ &    $-0.1365$  &  $0.1841$  &   $0.3167$   \\
 &Std error &$0.0349$&  $0.0347$  & $0.0334$  & $0.0320$  & $0.0439$  \\
 
 \hline

         $e^X$             &           Theory &                         &                  &         &                    0.2     &     0.4                     \\

         (rank $1$)  &    Simulation &  $-0.0486$  & $-0.0919$ &    $-0.0796$  &  $0.1432$  &   $0.2952$   \\
 &Std error &$0.0339$ &$0.0348$  & $0.0321$  & $0.0385$  & $0.0603$  \\
 \hline

         $I(X \leq 0.1)$             &           Theory &                         &                  &         &                    0.2     &     0.4                     \\

         (rank $1$)  &    Simulation &  $-0.1408$  & $-0.1342$ &    $-0.0961$  &  $0.1579$  &   $0.3124$   \\
 &Std error &$0.0316$ &$0.0319$  & $0.0326$  & $0.0325$  & $0.0371$  \\
\hline\hline
\end{tabular}  
\end{table}

\begin{table}[H]
\center
\caption{Average estimated memory parameters of some transformations of $2,000$ simulated  FARIMA$(1,d,0)$ and FARIMA$(1,d,1)$  processes with $2,000$ observations in each process and $t(10)$ innovations (except the transformation $e^x$, for which we use Gaussian innovations since $\mathbb{E}e^{2X}<\infty$ is required). The FARIMA$(1,d,0)$ processes have $\phi_1=-0.3$ and  the FARIMA$(1,d,1)$ processes have $\phi_1=-0.4,\; \theta_1=0.7$.}    
\label{tab:2}
\bigskip
\scriptsize
\begin{tabular}{ll  |  ll | ll}
\hline  \hline                                                     
  $K(X)$  and its   &          &\multicolumn{ 4}{c}{Memory parameter of the original series $X$}        \\
 power rank& &\multicolumn{2}{c}{FARIMA$(1,d,0)$} & \multicolumn{2}{c}{FARIMA$(1,d,1)$}\\
 &    &             $d=0.2$ &     $d=0.4$   &   $d=0.2$ &     $d=0.4$   \\
\hline

 $X$             &           Theory &                          0.2    &    0.4         &  0.2    &    0.4           \\

         (rank $1$)  &    Simulation &    $0.1624$  &   $0.3663 $ &  $0.2136$  &   $0.4188 $   \\
         & Std error & $0.0325$ & $0.0332$& $0.0329$  &$0.0317$  \\

\hline

         $X^2$             &           Theory &   0            &     0.3            &   0            &     0.3         \\

         (rank $2$)  &    Simulation &  $0.0212$  &   $0.2107 $& $0.0646$  &   $0.3007 $   \\
         & Std error & $0.0329$ & $0.0635$&$0.0376$  & $0.0615$  \\

\hline

         $X^3$             &           Theory &                         0.2     &     0.4     &    0.2     &     0.4    \\

         (rank $1$)  &    Simulation &   $0.0640$  &  $0.2173$ &    $0.1237$  &  $0.3083$  \\
 & Std error & $0.0350$ & $0.0542$ &$0.0374$  & $0.0563$  \\

\hline

         $X^4$             &           Theory &              0                     &     0.3        &              0                     &     0.3              \\

         (rank $2$)  &    Simulation &   $0.0084$  &  $0.1119$ &     $0.0356$  &  $0.2214$   \\
          &Std error & $0.0308$  & $0.0696$ & $0.0372$  & $0.0778$ \\

         \hline

         $X^3-3X$             &           Theory &      &     0.2         &      &     0.2             \\

         (rank $3$)  &    Simulation &     $0.0087 $  &   $0.2049 $  &   $0.0406$  &   $0.2582 $   \\
 &Std error & $0.0347$  & $0.0800$  & $0.0384$  & $0.0695$\\
\hline

         $X^4-6X^2$             &           Theory &                   &     0.1        &                   &     0.1                \\

         (rank $4$)  &    Simulation &    $0.0004 $  &$0.0425$ &   $0.0144 $ &   $0.1800 $  \\
          &Std error & $0.0302$  & $0.0633$   & $0.0378$  & $0.0905$ \\

\hline

         $\sin X$             &           Theory &                           0.2     &     0.4          &                           0.2     &     0.4             \\

         (rank $1$)  &    Simulation &    $0.1417$  &   $0.2996$   & $0.1868$ &   $0.2896$ \\
 &Std error & $0.0323$ & $0.0433$  & $0.0328$  & $0.0412$\\
 \hline
 $e^X$             &           Theory &                            0.2    &     0.4     &                            0.2    &     0.4                   \\

         (rank $1$)  &    Simulation &    $0.1042$  &  $0.2672$   &  $0.1482$ &   $0.2856$ \\
 &Std error & $0.0364$  & $0.0504$& $0.0392$  & $0.0696$  \\

 \hline
 $I(X \leq 0.1)$                &           Theory &                            0.2    &     0.4     &                            0.2    &     0.4                   \\

         (rank $1$)  &    Simulation &    $0.1145$  &  $0.2728$   &  $0.1712$ &   $0.3228$ \\
 &Std error & $0.0327$  & $0.0366$& $0.0331$  & $0.0374$  \\

\hline\hline
\end{tabular}  
\end{table}

The simulation study with these polynomial or non-polynomial transformations clearly confirms the theoretical results in  Theorem \ref{cor2} and Theorem \ref{thm3} for FARIMA$(0,d,0)$ processes with $-1<d<1/2$ or in general FARIMA$(p,d,q)$ processes with $0<d<1/2$.  

One can also compare the result in Table \ref{tab:1} with the simulation study performed in Dittmann and Granger (2002). They obtained the theoretical results for the memory parameters of the polynomial transformations of stationary Gaussian FARIMA$(0,d,0)$ processes. But there were no  theoretical results for the FARIMA$(p,d,q)$ processes with $p$ or $q$ not zero or for the non-polynomial transformations, such as $K(x)=\sin x$, $e^x$, even in the Gaussian case. They performed  simulation study for all the transformations in the Table \ref{tab:1} of Gaussian FARIMA$(0,d,0)$ processes.  
As expected, due to the heavy tail innovation,  the result in Table \ref{tab:1} is slightly worse than the one in Dittmann and Granger (2002). The innovation $t(10)$ used here has heavier tail than Gaussian innovation.
\section{Transformations of non-stationary processes}
In this section, we explore the memory properties of polynomial transformations of one type non-stationary processes. 
In the case  $1/2< d<3/2$,  a non-stationary process $X_n$ can be defined as 
the sum of a FARIMA$(0,d-1,0)$ processes, i.e.,
\begin{align}
&X_n=X_0+\sum_{j=1}^{n}Y_j, \label{X0}
\end{align}
where the distribution of the random variable $X_0$ does not depend on $n$, 
\begin{align}
& Y_j=\sum_{i=0}^{\infty}a_i \varepsilon_{j-i}, \label{Y0}
\end{align}
$a_i=\frac{\Gamma(i+d-1)}{\Gamma(i+1) \Gamma(d-1)}$
and $\varepsilon_t$ are i.i.d. random variables with mean $0$ and variance $1$. $a_i\sim i^{d-2}/\Gamma(d-1)$ for large $i\in \mathbb{N}$.
As in Velasco (1999 a, b), one can define $X_n$ analogously in the case $d\ge 3/2$.   $X_n$ defined in this way is called Type I process, see e.g., Shao and Wu (2007). 

In the following theorem we obtain the memory property of $X_n^2$ for Type I processes $X_n$. The memory property of $K(X_n)$ with a general transformation $K(\cdot)$ of  Type I processes $X_n$ is complicate and we leave it as an open question.

Notice that  $X_n-X_{n-1}=Y_{n}$ is a FARIMA$(0,d-1,0)$ process.  Thus $X_n \sim LM(d)$.
In the following theorem,  we show that this is also true asymptotically for $X_n^2$ in the case $1/2< d<1$.
\begin{theorem}\label{thm4}
Let $X_n$ be a Type I non-stationary process with $1/2< d<1$.  Assume that $X_0=0$ and $\mathbb{E}\varepsilon^4<\infty$. Then $X_n^2$ is asymptotically LM($d$).
\end{theorem}
Theorem 3.1 shows that taking the square of a non-stationary long memory process does not change the size of the long memory parameter, which is contrast to the result of stationary FARIMA$(p,d,q)$ processes.  

The simulation study in Table \ref{tab:4} is to confirm the result in Theorem \ref{thm4}. 
We simulate FARIMA$(0,d-1,0)$ processes for each of the $d$ values, $d=0.55$, $0.65$, $0.75$, $0.85$, $0.95$. The i.i.d. innovations have Student $t$ distribution with degree of freedom $5$. By the Definition (\ref{X0}), the partial sum gives a Type I process. The method to produce FARIMA$(0,d-1,0)$ processes and the method to estimate the memory parameters are same as the ones in Section \ref{main}.  It is clear from Table \ref{tab:4} that $X_t^2$ is asymptotically LM$(d)$ process. 
Also the rest part of Table \ref{tab:4} seems to confirm one conjecture: under suitable moment conditions, any polynomial transformations of non-stationary FARIMA$(0,d,0)$ processes are LM$(d)$ processes, $1/2<d<1$. The memory property of polynomial transformations of FARIMA$(0,d,0)$ processes is not related to the power ranks of these transformations. 
\begin{table}[H]
\center
\caption{Average estimated parameters of some polynomial transformations of $2,000$ simulated FARIMA$(0,d,0)$ processes with $2,000$ observations in each process. For the transformation $K(x)=x^2$, the innovation of the original process $X_n$ has Student $t$ distribution with degree of freedom $5$. For other transformations, the innovation of the original process $X_n$ has Student $t$ distribution with degree of freedom $10$. }
\label{tab:4}
\bigskip
\scriptsize
\begin{tabular}{ll | lllll}
\hline\hline                                                      
$K(X)$  &  &  \multicolumn{ 5}{c}{Memory parameter $d$ of the original series $X$}     \\  
             &  &  $0.55$ &   $0.65$ &  $0.75$ &  $0.85$ & $0.95$  \\
\hline

$X^2$   &  Theory & 0.55 &  0.65    & 0.75 & 0.85 & 0.95    \\

&  Simulation & 0.4826 &  0.6170  & 0.7433 & 0.8629 &  $0.9647 $ \\
&Std error & $0.0673$  & $0.0592$  & $0.0530$  & $0.0531$ &   $0.0456$ \\

\hline

$X^3$  &  Simulation & 0.4845 &  0.6102  & 0.7308 & 0.8479 &  $0.9598 $ \\
&Std  error & $0.0637$  & $0.0639$  & $0.0620$  & $0.0597$ &   $0.0554$ \\

\hline

$X^4$ &  Simulation & 0.4265 &  0.5702  & 0.7066 & 0.8087 &  $0.9441 $ \\
&Std  error & $0.0913$  & $0.0873$  & $0.0821$  & $0.0757$ &   $0.0678$ \\

\hline

$X^3-3X$ &  Simulation & 0.4667 &  0.6029  & 0.7298 & 0.8502 &  $0.9564 $ \\
&Std  error & $0.0781$  & $0.0671$  & $0.0637$  & $0.0628$ &   $0.0554$ \\

\hline

$X^4-6X^2$ &  Simulation & 0.4037 &  0.5585  & 0.7024& 0.8280 &  $0.9452 $ \\
&Std  error & $0.1053$  & $0.0969$  & $0.0854$  & $0.0785$ &   $0.0675$ \\

\hline\hline
\end{tabular}  
\end{table}

\begin{remark}
The decomposition in Section \ref{main} is not applicable for the non-stationary process with $1/2<d<1$. In fact, the result in Theorem \ref{thm4} is not related to the ranks of the transformations.
\end{remark}

\section{Application in option processes}
The transformation $K(x)=(x-C)^+$ itself has independent interest. It is $x-C$ if $x\ge C>0$. Otherwise it is $0$. Notice that this $K(x)$ is not differentiable at $C$.  For the reason to be clear later, let $X\ge 0$ be a random variable with mean $\mu$. Then $Y=X-\mu$ has mean $0$ and \begin{align*}
K(X)=(X-C)^+=(Y-(C-\mu))^+:=H(Y).
\end{align*}
Let $G(y)$ be the distribution function of $Y$. Assume that the density function of $Y$ exists and let it be $g(y)$. Then 
\begin{align*}
H_\infty(y)&=\int (y+z-(C-\mu))^+dG(z)\\
&=\int _{C-\mu-y}^\infty(y+z-(C-\mu))dG(z)\\
&=(y-C+\mu)[1-G(C-\mu-y)]+\int _{C-\mu-y}^\infty zdG(z)
\end{align*}
and 
\begin{align*}
&H'_\infty(y)\\
&=1-G(C-\mu-y)+(y-C+\mu)g(C-\mu-y)+(C-\mu-y)g(C-\mu-y)\\
&=1-G(C-\mu-y).
\end{align*}

We then have $H'_\infty(0)=1-G(C-\mu)\ne 0$ if $C-\mu$ is small enough. Therefore in this case the power rank of this $H(\cdot)$ is $1$.  If a larger $C-\mu>0$ is in the range such that $H'_\infty(0)=1-G(C-\mu)\approx 0$ and at the same time $g(C-\mu)$ is quite far away from $0$, we can say that the power rank of this $H(\cdot)$ is $2$ since $H''_\infty(y)=g(C-\mu-y)$ and therefore $H''_\infty(0)=g(C-\mu)>0$. If $G(y)$ is smooth enough (measured by the order of differentiability), we can continue this procedure to find the power rank of $H(\cdot)$ for $C-\mu$ in different ranges. In particular, if $G^{(r)}(y)$ exists for any $r\in \mathbb{N}$, and for $C-\mu>0$ large enough, $G^{(r)}(C-\mu)\approx 0$, then we say the power rank is $\infty$ since $H^{(k)}_\infty(y)=(-1)^k g^{(k-2)}(C-\mu-y)$ and therefore $H^{(k)}_\infty(0)=(-1)^k g^{(k-2)}(C-\mu)\approx 0$. 

We conduct simulation study for transformations $(X_n-C)^+$ with $C=0.3$, $1.5$, $5$, $9$, $44.8$ and $45.5$ where $Y_n=X_n-\mu$ are FARIMA$(0,d,0)$ processes with $d=0.2$ and $0.4$ as in Section \ref{main}. The innovations of $X_n$ are the absolute values of Student $t$ random variables with degree of freedom $10$.  The way to estimate the memory parameters of the transformations including the selection of the bandwidth $[n^{4/5}]$ is same as the one in Section \ref{main}. 
For each $d=0.2$ and $0.4$, we conduct $N=2,000$ simulations with $n=2^{20}$ observations in each process. Notice that the mean $\mu$ of $X_n$ changes for different memory parameters $d$. Therefore the power rank of $(X_n-C)^+=(Y_n-(C-\mu))^+$ also varies with $\mu$ for each fixed $C>0$. 
The result is listed in Table \ref{tab:5}. Again, If $d=0.2$,  the theoretical memory parameters of transformations  of the FARIMA$(0,d,0)$ processes with rank greater than 2 are left in blank. In the table, there are no estimates if $C-\mu>5$ since the length of each simulated process $X_n$ is finite and therefore the transformed values are all zeros if $C-\mu$ is too large. We use NA to denote them. Theoretically, the memory of a degenerate time series is zero. The simulation study confirms the results in Theorem \ref{cor2} and the above analysis.

The study of the memory parameter of $K(x)=(x-C)^+$ has direct application to call option time series in finance. Suppose $X_n$ is the price process of the underlying asset and $C$ is the strike price, then $K(X_n)$ is the value of the call option. Our result shows that the memory parameter of $(X_n-C)^+$ is same as the memory parameter of the underlying asset $X_n$ if $C-\mu$ is small. The power rank of 
$(X_n-C)^+=(Y_n-(C-\mu))^+$ is $2$ approximately if $C-\mu$ is in some moderate range. In this case, according to  Theorem \ref{cor2} (which is confirmed by the simulation study), the memory parameter of $K(X_n)$ is $2d-1/2$ if the memory parameter $d$ of the original mean adjusted asset price process $X_n-\mu\sim$ FARIMA$(0,d,0)$ satisfies $1/4<d<1/2$.

Similar analysis can be conducted for the truncation function $K(x)=(C-x)^+$ and the put option time series $(C-X_n)^+$ at different $C>0$. 
\begin{table}[H]
\center
\caption{Average estimated memory parameter of some transformations $(X_n-C)^+$ of $2,000$ simulated  FARIMA$(0,d,0)$ processes $X_n-\mu$ with $|t(10)|$ innovations and $2^{20}$ observations in each process $X_n$.  $C=0.3$, $1.5$, $5$, $9$, $44.8$ and $45.5$ and  $Y_n=X_n-\mu$.}
\label{tab:5}
\bigskip
\scriptsize
\begin{tabular}{l  |  lllll}
\hline  \hline                                                     
      &\multicolumn{ 4}{c}{Memory parameter of the original series $X$}        \\
     &            $d=0.2$ & &    $d=0.4$         \\
 \hline

       $(Y-(0.3-\mu))^+$ &$(Y+6.89)^+$ & & $(Y+44.74)^+$\\
Rank &  $1$ & &$1$\\
                           Theory &                      0.2     &  &   0.4                     \\

         Simulation &        $0.1999$  &  &$0.3998$\\
 Std error &   $0.0024$ &&$0.0025$  \\

 \hline

       $(Y-(1.5-\mu))^+$ &$(Y+5.69)^+$ &&$(Y+43.24)^+$\\
 Rank&  $1$ && $1$\\
                          Theory &                        0.2     &  &   0.4     \\

              Simulation &   $0.1999$  &  & $0.3998$   \\
 Std error & $0.0025$  & &$0.0024$  \\

 \hline

        $(Y-(5-\mu))^+$ &$(Y+2.19)^+$ & &$(Y+39.74)^+$\\
Rank& $1$ & &$1$\\
                         Theory &                                      0.2     &   &  0.4                   \\

              Simulation &   $0.1999$  &  & $0.3999$   \\
 Std error & $0.0025$  & &$0.0024$  \\
 \hline

        $(Y-(9-\mu))^+$ &$(Y-1.81)^+$ & &$(Y+35.74)^+$\\
Rank&  $2$ & &$1$\\
                           Theory &          0             &  &   0.4                   \\

             Simulation &  $0.0246$  &   &$0.3998$   \\
 Std error &$0.0027$  & &$0.0025$  \\

 \hline

         $(Y-(44.8-\mu))^+$ &$(Y-37.61)^+$ & &$(Y-0.06)^+$\\
Rank&  $\infty$ & &$1$\\
                            Theory &                                &   &  0.4                   \\

               Simulation &   $NA$  &   &$0.3332$   \\
 Std error & $NA$  & &$0.0049$  \\
 \hline

        $(Y-(45.5-\mu))^+$ &$(Y-38.31)^+$ & &$(Y-0.76)^+$\\
Rank& $\infty$ & &$2$\\
                           Theory &                   &             &     0.3                   \\

              Simulation &   $NA$  &  & $0.2541$   \\
 Std error&   $NA$  &   &$0.0074$   \\
%
%
%
\hline\hline
\end{tabular}  
\end{table}

\section{Proofs}
\label{proof}
Define the projection operator $$\mathcal{P}_{i}X=\mathbb{E}(X|\mathcal{F}_i)-\mathbb{E}(X|\mathcal{F}_{i-1}).$$ 
We adopt the notations from Wu (2006) as follows: For  $j \geq 2,$  let $A_n(j)=\sum_{t=n}^{\infty}|a_t|^j,$ $\theta_n=|a_{n-1}|[|a_{n-1}|+A^{1/2}_n(4)+A^{k/2}_{n}(2)].$  

Recall that 
$K(\cdot)$ is a measurable function with the power rank of $k$.
By  Ho and Hsing (1997), $K(X_n)-\mathbb{E}K(X_n)$ can be deccomposed as  $U(\mathcal{F}_n)+S(\mathcal{F}_n)$, where
\begin{align*}
U(\mathcal{F}_n)= K_{\infty}^{(k)}(0)\sum_{0\leq j_1<j_2<\cdots<j_k<\infty} \prod_{s=1}^{k}a_{j_s} \varepsilon_{n-j_s}
\end{align*}
and
\begin{align*}
S(\mathcal{F}_n)=K(X_n)-\mathbb{E}K(X_n)-U(\mathcal{F}_n).
\end{align*}
We also have the  decomposition $ U(\mathcal{F}_{n+h})+S(\mathcal{F}_{n+h})$ for $K(X_{n+h})$. 
Therefore, $\mbox{Cov} \left(K(X_n),K(X_{n+h})\right)$ can be represented by 
\begin{align}\label{cov}
&\mbox{Cov}(U(\mathcal{F}_n), U(\mathcal{F}_{n+h}))+\mbox{Cov}(U(\mathcal{F}_n), S(\mathcal{F}_{n+h}))\notag\\
&+\mbox{Cov}( S(\mathcal{F}_{n}), U(\mathcal{F}_{n+h}))+\mbox{Cov}( S(\mathcal{F}_{n}), S(\mathcal{F}_{n+h})).
\end{align}

\subsection{Proofs of Theorem \ref{thm1}, Corollary \ref{corollary1} and Theorem \ref{cor2}}
We first find the bounds of $\mbox{Cov}( U(\mathcal{F}_{n}), U(\mathcal{F}_{n+h}))$, which are useful in the proofs of Theorem \ref{thm1}, Corollary \ref{corollary1} and Theorem \ref{cor2}.  
By Stirling's approximation, $\left( \begin{array}{c}      
h\\
k
\end{array} \right) \simeq h^{k}/k!$,  if $h$ is large enough. If $j_k\le h$, the quantity $\prod_{s=1}^{k}L(j_s)j_s^{-\beta}L(h+j_s)(h+j_s)^{-\beta}$ at least has order $h^{-2k\beta}L^k(h)\min_{1\le s\le h}L^k(s)$.
Therefore for the lower bound, 
\begin{align}\label{UU}
&\mbox{Cov}( U(\mathcal{F}_{n}), U(\mathcal{F}_{n+h}))= [K_{\infty}^{(k)}(0)]^2 \sum_{0\leq j_1<j_2<\cdots<j_k<\infty} \prod_{s=1}^{k}a_{j_s} a_{h+j_s}\\
&\ge [ K_{\infty}^{(k)}(0)]^2 \sum_{1\leq j_1<j_2<\cdots<j_k\le h} \prod_{s=1}^{k}L(j_s)j_s^{-\beta}L(h+j_s)(h+j_s)^{-\beta}\nonumber \\
&\gg [ K_{\infty}^{(k)}(0)]^2 \left( \begin{array}{c}
h\\
k
\end{array} \right)  h^{-2k\beta}L^k(h)\min_{1\le s\le h}L^k(s) \nonumber \\
&\simeq (k!)^{-1} [ K_{\infty}^{(k)}(0) ]^2 h^{k}h^{-2k\beta}L^k(h)\min_{1\le s\le h}L^k(s) \nonumber \\
&= (k!)^{-1}[ K_{\infty}^{(k)}(0)]^2 h^{k(1-2\beta)}L^k(h)\min_{1\le s\le h}L^k(s).\nonumber 
\end{align}
On the other hand,
\begin{align}\label{xi3}
&\mbox{Cov}( U(\mathcal{F}_{n}), U(\mathcal{F}_{n+h})) \\
&\leq [ K_{\infty}^{(k)}(0)]^2 (\sum_{0\le i< \infty}a_ia_{h+i})^k \nonumber \\
&= [K_{\infty}^{(k)}(0)]^2(\sum_{i=0}^{h}a_ia_{h+i}+\sum_{i=h+1}^{\infty}a_ia_{h+i})^k \nonumber \\
&=\mathcal{O}[h^{k(1-2\beta)}L^{2k}(h)].\nonumber
\end{align}

\noindent {\bf Proof of Theorem \ref{thm1}}\\
\noindent  We shall  estimate  the other covariances in (\ref{cov}) for each case. In 
the case that $k(2\beta-1)<1$, we first apply the projection operator to the terms in the covariances and then apply the Cauchy-Schwarz inequality. Then
\begin{align}\label{SS}
&|\mbox{Cov}( S(\mathcal{F}_{n}),S(\mathcal{F}_{n+h}))| \nonumber\\
&=|\sum_{i=-\infty}^{n-1}\sum_{j=-\infty}^{n+h-1}\mbox{Cov}(\mathcal{P}_{i+1} S(\mathcal{F}_{n}),\mathcal{P}_{j+1} S(\mathcal{F}_{n+h}))|\nonumber \\
&=|\sum_{i=-\infty}^{n-1}\sum_{j=-\infty}^{n+h-1}\mathbb{E}[\mathcal{P}_{i+1} S(\mathcal{F}_{n})\mathcal{P}_{j+1} S(\mathcal{F}_{n+h})]|\nonumber \\
&=|\sum_{i=-\infty}^{n-1}\mathbb{E}[\mathcal{P}_{i+1} S(\mathcal{F}_{n})\mathcal{P}_{i+1} S(\mathcal{F}_{n+h})]| \\  \label{EE}
&\leq \sum_{i=-\infty}^{n-1} \|\mathcal{P}_{i+1} S(\mathcal{F}_{n})\| \|\mathcal{P}_{i+1} S(\mathcal{F}_{n+h})\|\nonumber \\
&=\sum_{i=-\infty}^{n-1} \|\mathcal{P}_{1} S(\mathcal{F}_{n-i})\| \|\mathcal{P}_{1} S(\mathcal{F}_{n+h-i})\|\nonumber \\
&=\sum_{i=-\infty}^{n-1} \mathcal{O} (\theta_{n-i} \theta_{n+h-i})=\sum_{i=1}^{\infty} \mathcal{O} (\theta_{i} \theta_{i+h}).
\end{align}
Equality  (\ref{SS}) is true because if $i \neq j,$ suppose $i <j,$ then 
\begin{align*}
&\mathbb{E}[\mathcal{P}_{i+1} S(\mathcal{F}_{n})\mathcal{P}_{j+1} S(\mathcal{F}_{n+h})] =\mathbb{E}\big \{\mathbb{E}[\mathcal{P}_{i+1} S(\mathcal{F}_{n})\mathcal{P}_{j+1} S(\mathcal{F}_{n+h})|\mathcal{F}_{i+1}] \big \}\nonumber \\
&=\mathbb{E}\big \{\mathcal{P}_{i+1} S(\mathcal{F}_{n}) \mathbb{E}[\mathcal{P}_{j+1} S(\mathcal{F}_{n+h})|\mathcal{F}_{i+1}] \big \}\nonumber \\
&=\mathbb{E}\big \{\mathcal{P}_{i+1} S(\mathcal{F}_{n}) \mathbb{E} \big [(\mathbb{E}[S(\mathcal{F}_{n+h})|\mathcal{F}_{j+1}]-\mathbb{E}[S(\mathcal{F}_{n+h})|\mathcal{F}_{j}])|\mathcal{F}_{i+1}\big ] \big \}\nonumber \\
&=\mathbb{E}\big \{\mathcal{P}_{i+1} S(\mathcal{F}_{n}) \big ( \mathbb{E}[S(\mathcal{F}_{n+h})|\mathcal{F}_{i+1}]-\mathbb{E}[S(\mathcal{F}_{n+h})|\mathcal{F}_{i+1}]\big ) \big \}=0.\nonumber 
\end{align*}
Equality (\ref{EE}) is the result of Theorem 5 (Reduction principle) of Wu (2006). By Karamata's theorem  (Seneta, 1976), $A_{n}(j)=\mathcal{O}[n^{1-j\beta}L^{j}(n)]$ for $j \geq 2.$ Therefore, under the condition $k(2\beta-1)<1$,
\begin{align*}
&\sum_{i=1}^{\infty} \mathcal{O} (\theta_{i} \theta_{i+h})\nonumber \\
&=\sum_{i=0}^{\infty}a_i a_{i+h} \mathcal{O}\big \{[(i+1)^{1-2\beta} L^{2}(i+1)(i+h+1)^{1-2\beta} L^{2}(i+h+1)]^{k/2}\big \}\;\;\nonumber \\
&=\sum_{i=1}^{\infty}\mathcal{O}[i^{-\beta+k(1-2\beta)/2}(i+h)^{-\beta+k(1-2\beta)/2} L^{k+1}(i) L^{k+1}(i+h)]\nonumber\\
&=\sum_{i=1}^{h}\mathcal{O}[i^{-\beta+k(1-2\beta)/2}h^{-\beta+k(1-2\beta)/2}(1+i/h)^{-\beta+k(1-2\beta)/2} L^{k+1}(i) L^{k+1}(i+h)]\nonumber\\
&+\sum_{i=h+1}^{\infty}\mathcal{O}[i^{-2\beta+k(1-2\beta)}(1+h/i)^{-\beta+k(1-2\beta)/2} L^{k+1}(i) L^{k+1}(i+h)].
\end{align*}
Applying Karamata's theorem again, we have 

(i) $\sum_{i=1}^{\infty} \mathcal{O} (\theta_{i} \theta_{i+h})=\mathcal{O}[L^{k+1}(h) h^{\frac{k(1-2\beta)-2\beta}{2}}]$\;\;\;\; if $(k+1)(2\beta-1)>1;$

(ii)$\sum_{i=1}^{\infty} \mathcal{O} (\theta_{i} \theta_{i+h})=\mathcal{O}[L^{2k+2}(h) h^{(k+1)(1-2\beta)}]$\;\; if $(k+1)(2\beta-1)< 1$;

(iii) $\sum_{i=1}^{\infty} \mathcal{O} (\theta_{i} \theta_{i+h})=\max \big \{\mathcal{O}[L^{2k+2}(h) h^{-1}], \mathcal{O}[L^{k+1}(h) h^{-1}\sum_{i=1}^h i^{-1}L^{k+1}(i)]\big \}$ \;\;\;
 if $(k+1)(2\beta-1)= 1$.

\noindent By the calculation in (\ref{UU}), each of the above terms (i), (ii) and (iii) is less than  $ \mbox{Cov}( U(\mathcal{F}_{n}), U(\mathcal{F}_{n+h}))$. Thus,
\begin{align*}
|\mbox{Cov}( S(\mathcal{F}_{n}),S(\mathcal{F}_{n+h}))|<\mbox{Cov}( U(\mathcal{F}_{n}),U(\mathcal{F}_{n+h})).
\end{align*}

With the same arguments as in $|\mbox{Cov}( S(\mathcal{F}_{n}),S(\mathcal{F}_{n+h}))|$, 
\begin{align}
|\mbox{Cov}(U(\mathcal{F}_{n}),S(\mathcal{F}_{n+h}))|
=\sum_{i=1}^{\infty}  \|\mathcal{P}_{1} U(\mathcal{F}_{i})\| \mathcal{O} ( \theta_{i+h}).\nonumber
\end{align}
Obvioulsly, $\mathcal{O} ( \theta_{i+h})=\mathcal{O}[(i+h)^{-\beta+\frac{k(1-2\beta)}{2}}L^{k+1}(h+i)].$ Also
\begin{align*}
&\mathcal{P}_{1} U(\mathcal{F}_{i})/K_\infty^{(k)}(0)\nonumber \\
&=\mathbb{E}\big (\sum_{0\leq j_1< \cdots < j_k<\infty} \prod_{s=1}^k a_{j_s}\varepsilon_{i-j_s}|\mathcal{F}_1 \big)-\mathbb{E}\big (\sum_{0\leq j_1< \cdots < j_k<\infty} \prod_{s=1}^k a_{j_s}\varepsilon_{i-j_s}|\mathcal{F}_0\big )\nonumber \\
&=\sum_{i-1\leq j_1< \cdots < j_k<\infty} \prod_{s=1}^k a_{j_s}\varepsilon_{i-j_s}-\sum_{i\leq j_1< \cdots < j_k<\infty} \prod_{s=1}^k a_{j_s}\varepsilon_{i-j_s}\nonumber \\
&=a_{i-1}\varepsilon_{1}\sum_{i\leq j_2< \cdots < j_k<\infty} \prod_{s=2}^k a_{j_s}\varepsilon_{i-j_s}. 
\end{align*}
Then 
\begin{align*}
&||\mathcal{P}_{1} U(\mathcal{F}_{i})||^2=\mathbb{E}[\mathcal{P}_{1} U(\mathcal{F}_{i})]^2\\
&=[K_\infty^{(k)}(0)]^2a^{2}_{i-1}\mathbb{E}(\varepsilon_1^2)\mathbb{E}\big (\sum_{i\leq j_2< \cdots < j_k<\infty} \prod_{s=2}^k a_{j_s}\varepsilon_{i-j_s}\big )^2\nonumber \\
&\leq [K_\infty^{(k)}(0)]^2a^{2}_{i-1}A^{k-1}_i(2)[\mathbb{E}(\varepsilon_1^2)]^{k}.
\end{align*}
See also Wu (2006). Hence
\begin{align*}
||\mathcal{P}_{1} U(\mathcal{F}_{i})||=\mathcal{O}[i^{-\beta+(k-1)(1-2\beta)/2}L^{k}(i)].\nonumber \\
\end{align*}
In consequence, using Karamata's  theorem,  we have 
\begin{align*}
&\sum_{i=1}^{\infty}  ||\mathcal{P}_{1} U(\mathcal{F}_{i})|| \mathcal{O} ( \theta_{i+h})\nonumber\\ &=\sum_{i=1}^{\infty}\mathcal{O}\big [(i+h)^{-\beta+\frac{k(1-2\beta)}{2}}L^{k+1}(h+i) i^{-\beta+(k-1)(1-2\beta)/2}L^{k}(i)\big]\nonumber \\
&=h^{-\beta+\frac{k(1-2\beta)}{2}}\sum_{i=1}^{h}\mathcal{O}\big [(1+i/h)^{-\beta+\frac{k(1-2\beta)}{2}}L^{k+1}(h+i) i^{-\beta+(k-1)(1-2\beta)/2}L^{k}(i)\big ]\nonumber \\
&+\sum_{i=h+1}^{\infty}\mathcal{O}\big [ i^{-2\beta+(2k-1)(1-2\beta)/2}(1+h/i)^{-\beta+\frac{k(1-2\beta)}{2}}L^{k+1}(h+i)L^{k}(i)\big ]\nonumber \\
&=\mathcal{O}[h^{(1-2\beta)(k+1/2)}L^{2k+1}(h)],
\end{align*}
which is less than $ \mbox{Cov}( U(\mathcal{F}_{n}), U(\mathcal{F}_{n+h}))$.
So $$|\mbox{Cov}(U(\mathcal{F}_{n}),S(\mathcal{F}_{n+h}))|<\mbox{Cov}(U(\mathcal{F}_{n}),U(\mathcal{F}_{n+h})).$$

Similarly, 
$|\mbox{Cov}(U(\mathcal{F}_{n+h}),S(\mathcal{F}_{n}))|= \sum_{i=1}^{\infty}  \|\mathcal{P}_{1}U(\mathcal{F}_{i+h})\| \mathcal{O} ( \theta_{i}),$
which is 

(i)$\mathcal{O}\big [h^{(k+\frac{1}{2})(1-2\beta)}L^{2k+1}(h)\big]$ if $(k+1)(2\beta-1)<1$;

(ii) $\mathcal{O}\big[h^{-\beta+\frac{k-1}{2}(1-2\beta)}L^{k}(h)\big ]$ if $(k+1)(2\beta-1)>1;$

(iii) $\max\big \{\mathcal{O}[h^{-\frac{2k+1}{2(k+1)}}L^k(h)\sum_{i=1}^h i^{-1}L^{k+1}(i)], \mathcal{O}[h^{-\frac{2k+1}{2(k+1)}}L^{2k+1}(h)]\big \}$ if $(k+1)(2\beta-1)=1$. 

\noindent Each of the terms (i), (ii) and (iii) is less than $\mbox{Cov}(U(\mathcal{F}_{n}), U(\mathcal{F}_{n+h}))$ by the analysis in (\ref{UU}). Hence 
$|\mbox{Cov}(U(\mathcal{F}_{n+h}),S(\mathcal{F}_{n}))|<\mbox{Cov}(U(\mathcal{F}_{n}),U(\mathcal{F}_{n+h})).$

So under the condition $k<(2\beta-1)^{-1}$,  
$$\mbox{Cov}(K(X_n), K(X_{n+h})) \simeq  Cov(U(\mathcal{F}_{n}),U(\mathcal{F}_{n+h})).$$
 But by (\ref{UU}), 
$$\mbox{Cov}(U(\mathcal{F}_{n}),U(\mathcal{F}_{n+h})) \ge (k!)^{-1}[ K_{\infty}^{(k)}(0)]^2 h^{k(1-2\beta)}L^k(h)\min_{1\le s\le h}L^k(s),$$ which is not summable. 
Therefore $K(X_n)$ has long memory in the covariance sense if $k(2\beta-1)<1$.

Now, we consider the case $k > (2\beta-1)^{-1}$. With similar arguments as the case $k(2\beta-1)<1$, we have 
\begin{align}
&|\mbox{Cov}(U(\mathcal{F}_{n}),S(\mathcal{F}_{n+h}))|= \mathcal{O}[h^{-2\beta} L^{2}(h)]\;\;\; \text{if}\; (k-1)(2\beta-1)>1,\nonumber\\
&|\mbox{Cov}(U(\mathcal{F}_{n}),S(\mathcal{F}_{n+h}))|= \mathcal{O}\big \{h^{-2\beta} \max [L^{2}(h), L^{k+1}(h)]\big\}\;\;\; \text{if}\; (k-1)(2\beta-1)=1,\nonumber\\
&|\mbox{Cov}(U(\mathcal{F}_{n}),S(\mathcal{F}_{n+h}))|= \mathcal{O}[h^{-\beta+k(1-2\beta)/2} L^{k+1}(h)]\; \;\;\text{if} \;(k-1)(2\beta-1)<1;\label{us}
\end{align}
\begin{align}
&|\mbox{Cov}(S(\mathcal{F}_{n}),S(\mathcal{F}_{n+h}))|= \mathcal{O}[h^{-2\beta} L^{2}(h)]\;\;\; \text{if}\; (k-1)(2\beta-1)>1,\nonumber\\
&|\mbox{Cov}(S(\mathcal{F}_{n}),S(\mathcal{F}_{n+h}))|= \mathcal{O}\{h^{-2\beta} \max[L^{2}(h), L^{k+1}(h)]\}\;\;\; \text{if}\; (k-1)(2\beta-1)=1,\nonumber\\
&|\mbox{Cov}(S(\mathcal{F}_{n}),S(\mathcal{F}_{n+h})|= \mathcal{O}[h^{-\beta+(1-2\beta)k/2} L^{k+1}(h)]\; \;\;\text{if} \;(k-1)(2\beta-1)<1; \label{ss}
\end{align}
and
\begin{align}
&|\mbox{Cov}(S(\mathcal{F}_{n}),U(\mathcal{F}_{n+h}))|= \mathcal{O}[h^{-\beta+(k-1)(1-2\beta)/2} L^{k}(h)],\label{su}
\end{align}
which are all summable. Additionally, 
(\ref{xi3}) is also summable in this case  $k(2\beta-1)>1$. Therefore $K(X_n)$ has short memory in the covariance sense if $k(2\beta-1)>1.$\\

\noindent {\bf Proof of Corollary \ref{corollary1}}\\
\noindent We just need consider the case that $(2\beta-1)^{-1}$ is an integer and the power rank $k=(2\beta-1)^{-1}$. In this case, since $\lim_{n\rightarrow \infty}L(n)=L$, by  (\ref{UU}), we have 
\begin{eqnarray}
\mbox{Cov}(U(\mathcal{F}_{n}),U(\mathcal{F}_{n+h}))\gg (k!)^{-1}[ K_{\infty}^{(k)}(0)]^2 h^{-1}.\nonumber 
\end{eqnarray}
Hence $K(X_n)$ has long memory in the covariance sense in this case. \\

\noindent {\bf Proof of Theorem \ref{cor2}}\\
\noindent Since  $\lim_{j\to\infty} \frac{a_j}{j^{d-1}}=\frac{\theta(1)}{\phi(1)\Gamma(d)}$, we have $\beta=1-d$ and the slowly varying function is a constant asymptotically. Hence, applying (\ref{UU}) and (\ref{xi3}) yields
\begin{align}\label{exactcov}
\mbox{Cov}(U(\mathcal{F}_{n}),U(\mathcal{F}_{n+h})) \simeq h^{k(1-2\beta)}
= h^{k(2d-1)} = h^{2\bar{d}-1}
\end{align}
with  $ \bar{d}=(d-1/2)k+1/2 > 0$.

We first consider the case  $k<(2\beta-1)^{-1}$. By (\ref{exactcov}) and the proof of  Theorem \ref{thm1},  
$\mbox{Cov}(K(X_n), K(X_{n+h})) \simeq  \mbox{Cov}(U(\mathcal{F}_{n}),U(\mathcal{F}_{n+h}))\simeq h^{k(1-2\beta)}.$   
Then by the same argument as in Proposition 1 of Dittmann and Granger (2002), $K(X_n)$ is a long-memory process LM$(\tilde{d})$ when $k(2\beta-1)<1$.

In the case that  $k(2\beta-1)>1$ and  $(k-1)(2\beta-1)<1,$ from   (\ref{exactcov}),  (\ref{us}),  (\ref{ss}) and (\ref{su}),
 $\mbox{Cov}(K(X_n), K(X_{n+h}))$ is dominated by $\mbox{Cov}(U(\mathcal{F}_{n}),U(\mathcal{F}_{n+h}))$. So,
 $$\mbox{Cov}(K(X_n), K(X_{n+h})) \simeq  h^{2\bar{d}-1}$$ with  $ \bar{d}=(d-1/2)k+1/2<0$.  Therefore the process $K(X_n)$ has the same autocorrelation delay pattern as an FARIMA$(0, \bar{d},0)$ process.  But we shall show that  it is a short-memory LM$(0)$ process.
Denote $f_K(\lambda)$ as the spectral density of $K(X_n)$. Since $\mbox{Cov}(K(X_n), K(X_{n+h}))$ is dominated by $\mbox{Cov}( U(\mathcal{F}_{n}), U(\mathcal{F}_{n+h}))$,  which is positive and summable,  then
\begin{align*}
0<f_{K}(0)=\mbox{Var}(K(X_n))+2\sum_{h=1}^\infty \mbox{Cov}(K(X_n), K(X_{n+h}))<\infty.
\end{align*}
Therefore $K(X_n)$ is a LM$(0)$ process.

\subsection{Proofs of Theorem \ref{thm2}, \ref{thm3}}
\noindent {\bf Proof of Theorem \ref{thm2}}

Again, using the projection operator and the Cauchy-Schwarz inequality, we have 
\begin{align}
&\mbox{Cov}(K(X_n), K(X_{n+h}))=\mbox{Cov}(\sum_{i=-\infty}^{n-1}\mathcal{P}_{i+1}K(X_n), \sum_{j=-\infty}^{n+h-1}\mathcal{P}_{j+1}K(X_{n+h}))\nonumber \\
&=\sum_{i=-\infty}^{n-1}\mathbb{E}[\mathcal{P}_{i+1}K(X_n)\mathcal{P}_{i+1}K(X_{n+h})]
\leq\sum_{i=-\infty}^{n-1}\|\mathcal{P}_{i+1}K(X_n)\| \|\mathcal{P}_{i+1}K(X_{n+h})\| \nonumber \\
&= \sum_{i=-\infty}^{n-1}\|\mathcal{P}_{1}K(X_{n-i})\| \|\mathcal{P}_{1}K(X_{n+h-i})\| \nonumber \\
&=\sum_{i=-\infty}^{n-1}\mathcal{O}(|a_{n-i-1}a_{n+h-i-1}|) \label{th2}\\
&=\sum_{i=0}^{\infty} \mathcal{O}(|a_{i}a_{i+h}|),\nonumber
\end{align}
where equality (\ref{th2}) is obtained from Wu (2006): $\|\mathcal{P}_{1}K(X_n)\|= \mathcal{O}(|a_{n-1}|),$ if condition (\ref{cond2}) holds.
Thus
\begin{align*}
\sum_{h=1}^{\infty}\mbox{Cov}(K(X_n), K(X_{n+h}))= \sum_{h=1}^{\infty}\sum_{i=0}^{\infty} \mathcal{O}(|a_{i}a_{i+h}|),
\end{align*}
which is finite.
This finishes the proof. \\

\noindent {\bf Proof of Theorem \ref{thm3}}\\
\noindent Denote $f_K(\lambda)$ as the spectral density of $K(X_n)=X^2_n,$ then 
 \begin{align}
&f_{K}(0)=\mbox{Var}(K(X_n))+2\sum_{h=1}^\infty \mbox{Cov}(K(X_n), K(X_{n+h}))\notag\\
&=\sum_{i=0}^{\infty} a^{4}_{i}\mbox{ Var}(\varepsilon^2_1)+4\sum_{0\leq i<j<\infty}a^{2}_{i}a^{2}_{j} +2\sum_{h=1}^{\infty}\sum_{i=0}^{\infty} a^2_{i} a^2_{h+i}\mbox{Var}(\varepsilon^2_1) \notag\\
&+8\sum_{h=1}^{\infty}\sum_{0\leq i<j<\infty} a_{i} a_{h+i}a_{j} a_{h+j} .\notag\\
&=\sum_{i=0}^{\infty} a^{4}_{i} \mbox{Var}(\varepsilon^2_1)+4\sum_{1\leq i<j<\infty}a^{2}_{i}a^{2}_{j} +2\sum_{h=1}^{\infty}\sum_{i=0}^{\infty} a^2_{i} a^2_{h+i}\mbox{Var}(\varepsilon^2_1) \notag\\
&+8\sum_{h=1}^{\infty}\sum_{1\leq i<j<\infty} a_{i} a_{h+i}a_{j} a_{h+j}+4\sum_{i=1}^{\infty}a_i^2+8\sum_{h=1}^{\infty} a_h\sum_{i=1}^{\infty}a_i a_{i+h}.  \notag
 \end{align}

We shall show that $f_{K}(0)>0.$ 
The condition $-1<d<0$ implies $a_i<0$ for all $i>0$. Therefore only the last term of the above decomposition of $f_{K}(0)$ is negative. 
To prove $f_{K}(0)>0$, it suffices to show that 
\begin{align}
&\sum_{i=1}^{\infty}a_i^2+2\sum_{h=1}^{\infty} a_h\sum_{i=1}^{\infty}a_i a_{i+h} \label{fk0}
\end{align}
is positive.
In fact,     
 \begin{align}
& (\ref{fk0})=\sum_{i=1}^{\infty}a_i^2+2\sum_{h=1}^{\infty} a_h\sum_{i=1}^{\infty}a_i a_{i+h}\notag\\
&=\sum_{i=1}^{\infty}\frac{\Gamma^2(i+d)}{\Gamma^2(d)\Gamma^2(i+1)}+2\sum_{h=1}^{\infty} \frac{\Gamma(h+d)}{\Gamma(d)\Gamma(h+1)}\sum_{i=1}^{\infty}\frac{\Gamma(i+d)\Gamma(i+h+d)}{\Gamma^2(d)\Gamma(i+1)\Gamma(i+h+1)} \notag\\
&=\sum_{i=1}^{\infty}\frac{\Gamma^2(i+d)}{\Gamma^2(d)\Gamma^2(i+1)}+2\sum_{h=1}^{\infty} \frac{\Gamma^2(h+d)}{\Gamma^2(d)\Gamma^2(h+1)} [F(d,h+d;h+1;1)-1] \label{F=0} \\
&=\sum_{h=1}^{\infty}\frac{\Gamma^2(h+d)}{\Gamma^2(d)\Gamma^2(h+1)} [2F(d,h+d;h+1;1)-1]\notag\\
&=\sum_{h=1}^{\infty}\frac{\Gamma^2(h+d)}{\Gamma^2(d)\Gamma^2(h+1)}\frac{2\Gamma(h+1)\Gamma(1-2d)}{\Gamma(h+1-d)\Gamma(1-d)}-\sum_{h=1}^{\infty}\frac{\Gamma^2(h+d)}{\Gamma^2(d)\Gamma^2(h+1)}  \label{G=0}\\
&=\frac{2\Gamma(1-2d)}{\Gamma^2(1-d)}[F(d,d;1-d;1)-1]-[F(d,d;1;1)-1]\notag\\
&=\frac{2\Gamma(1-2d)}{\Gamma^2(1-d)}[\frac{\Gamma(1-d)\Gamma(1-3d)}{\Gamma^2(1-2d)}-1]-\frac{\Gamma(1-2d)}{\Gamma^2(1-d)}+1\label{G1}\\
&=\frac{3\Gamma(-3d)\Gamma(-d)-d\Gamma^2(-d)\Gamma(-2d)-6\Gamma^2(-2d)}{-d\Gamma(-2d)\Gamma^2(-d)}\label{P=0}.
 \end{align}
The notation $F(a,b;c;z)$ from (\ref{F=0}) and thereafter is the hypergeometric series. (\ref{G=0}) and (\ref{G1}) are obtained by applying the Gauss's theorem for hypergeometric series (Gauss, 1866), see also page 2 of Bailey (1935).
The denominate of the last equation (\ref{P=0}) is positive since $-1<d<0$. Hence it suffices to prove that $$3\Gamma(-3d)\Gamma(-d)-d\Gamma^2(-d)\Gamma(-2d)-6\Gamma^2(-2d)>0.$$
Define $f(x)=3\Gamma(3x)\Gamma(x)+x\Gamma^2(x)\Gamma(2x)-6\Gamma^2(2x),$ $0<x<1$. The function $f(x)$ is continuous for $x>0$. 
Straight forward numerical calculation shows that  $f(x)>1/4>0$ for all $0<x<1$. Thus,  (\ref{fk0}) is positive.
Hence, $X^2_n$ is a LM$(0)$ process.

\subsection{Proof of Theorem \ref{thm4}}

By (\ref{X0}) and (\ref{Y0}), $X_n$ can be written in the form 
\begin{eqnarray} 
X_n=\sum_{j=0}^{\infty}b_{n}(j) \varepsilon_{n-j},\nonumber
\end{eqnarray}
where $b_{n}(j)=\sum_{i=0}^{j}a_i,$ if $0\le j \leq n$, 
and $b_{n}(j)=\sum_{i=j-n+1}^{j}a_i,$ if $j>n$. By convention, define $b_n(j)=0$ for $j<0$.  Let $Z_n=X_n+X_{n-1}$. Then $X^{2}_n-X^{2}_{n-1}=Y_n Z_n$ and 
\begin{align}\label{zn}
Z_n=\sum_{j=0}^{\infty}[b_n(j)+b_{n-1}(j-1)]\varepsilon_{n-j}.
\end{align}
i). Denote $\gamma_y(h)=\mbox{Cov}(Y_n, Y_{n+h})$ as the autocovariance function of the stationary process $Y_n$. We first show that, in the case $d<5/4$, 
\begin{align}
&\mbox{Cov}(Y_nZ_n, Y_{n+h}Z_{n+h})\nonumber \\
&= \gamma_y(h)\mbox{Cov}(Z_n, Z_{n+h})+\mbox{Cov}(Y_n, Z_{n+h})\mbox{Cov}(Z_n, Y_{n+h})+C(n,h) \label{ine}
\end{align}
as $n \to \infty$ for some constant $C(n,h)$ with uniform bound $0<C<\infty$. In fact, by (\ref{zn}), 
$$Y_nZ_n=\sum_{i=0}^\infty\sum_{j=0}^{\infty}a_i[b_n(j)+b_{n-1}(j-1)]\varepsilon_{n-i}\varepsilon_{n-j}$$
and by the change of variables, 
$$Y_{n+h}Z_{n+h}=\sum_{i=-h}^\infty\sum_{j=-h}^{\infty}a_{i+h}[b_{n+h}(j+h)+b_{n+h-1}(j+h-1)]\varepsilon_{n-i}\varepsilon_{n-j}.$$
Hence by the independence of the innovations $\varepsilon_i$, $i\in \mathbb{Z}$, 
\begin{align*}
&\mbox{Cov}(Y_nZ_n, Y_{n+h}Z_{n+h})\nonumber \\
&=\sum_{i=0}^{\infty}\sum_{j=0}^{\infty}\left\{a_i a_{h+i}[b_n(j)+b_{n-1}(j-1)]\right.\\
&\left.\times[b_{n+h}(h+j)+b_{n+h-1}(h+j-1)]Var(\varepsilon_{n-i} \varepsilon_{n-j})\right\} \nonumber \\
&+\sum_{i=0}^{\infty}\sum_{j=0}^{\infty}\left\{a_i a_{h+j}[b_n(j)+b_{n-1}(j-1)]\right.\\
&\left.\times[b_{n+h}(h+i)+b_{n+h-1}(h+i-1)]Var(\varepsilon_{n-i} \varepsilon_{n-j})\right\}. \nonumber 
\end{align*}
On the other hand, 
\begin{eqnarray}
&&\gamma_y(h)\mbox{Cov}(Z_n, Z_{n+h})=\mbox{Cov}(Y_n, Y_{n+h})\mbox{Cov}(Z_n, Z_{n+h})\nonumber \\
&&=\sum_{i=0}^{\infty}\sum_{j=0}^{\infty}a_i a_{h+i}[b_n(j)+b_{n-1}(j-1)][b_{n+h}(h+j)+b_{n+h-1}(h+j-1)]\nonumber 
\end{eqnarray}
and 
\begin{align*}
&\mbox{Cov}(Y_n, Z_{n+h})\mbox{Cov}(Z_n, Y_{n+h})\\
&=\sum_{i=0}^{\infty}\sum_{j=0}^{\infty}a_i a_{h+j}[b_n(j)+b_{n-1}(j-1)][b_{n+h}(h+i)+b_{n+h-1}(h+i-1)].
\end{align*}
Since $Var(\varepsilon_1 \varepsilon_2)=1$, we have 
\begin{align}
&\mbox{Cov}(Y_nZ_n, Y_{n+h}Z_{n+h})- \gamma_y(h)\mbox{Cov}(Z_n, Z_{n+h})-\mbox{Cov}(Y_n, Z_{n+h})\mbox{Cov}(Z_n, Y_{n+h})\nonumber \\
&=2\sum_{j=0}^{\infty}a_j a_{h+j}[b_n(j)+b_{n-1}(j-1)][b_{n+h}(h+j)+b_{n+h-1}(h+j-1)]\nonumber\\
&\times[Var(\varepsilon^2)-1]\nonumber\\
&=2\sum_{j=0}^{n}a_j a_{j+h}\left[\sum_{i=0}^{j}a_i+\sum_{i=0}^{j-1}a_i\right]\left[\sum_{i=0}^{h+j}a_i+\sum_{i=0}^{h+j-1}a_i\right][Var(\varepsilon^2)-1] \nonumber\\
&+2\sum_{j=n+1}^{\infty}a_j a_{j+h}\left[\sum_{i=j-n+1}^{j}a_i+\sum_{i=j-n+1}^{j-1}a_i\right]\left[\sum_{i=j-n+1}^{h+j}a_i+\sum_{i=j-n+1}^{h+j-1}a_i\right]\nonumber\\
&\times[Var(\varepsilon^2)-1].\label{a2}
\end{align}
In the above equations,  since $a_i \sim \frac{i^{d-2}}{\Gamma(d-1)}$,
\begin{align*}
\left|a_j a_{j+h}\left[\sum_{i=0}^{j}a_i+\sum_{i=0}^{j-1}a_i\right]\left[\sum_{i=0}^{h+j}a_i+\sum_{i=0}^{h+j-1}a_i\right]\right| <C j^{2d-4}
\end{align*}
for some $0<C<\infty$. 
Similarly,
\begin{align*}
\left|a_j a_{j+h}\left[\sum_{i=j-n+1}^{j}a_i+\sum_{i=j-n+1}^{j-1}a_i\right]\left[\sum_{i=j-n+1}^{h+j}a_i+\sum_{i=j-n+1}^{h+j-1}a_i\right]\right|<Cj^{4d-6}
\end{align*}
for some $0<C<\infty$.  So that (\ref{a2}) converges as $n \to \infty$ if $d<5/4.$ Consequently, in this case $d<5/4$, (\ref{ine}) holds as $n \to \infty$.

ii). Now we consider the second term of (\ref{ine}), $\mbox{Cov}(Y_n, Z_{n+h})\mbox{Cov}(Z_n, Y_{n+h})$.
\begin{align*}
&|\mbox{Cov}(Y_n, Z_{n+h})|=\left|\sum_{j=0}^{\infty}a_j[b_{n+h}(h+j)+b_{n+h-1}(h+j-1)]\right|\\
&=\left|\sum_{j=0}^{n}a_j[\sum_{i=0}^{h+j}a_i+\sum_{i=0}^{h+j-1}a_i] 
+\sum_{j=n+1}^{\infty}a_j[\sum_{i=j-n+1}^{h+j}a_i+\sum_{i=j-n+1}^{h+j-1}a_i]\right| \\
&<C\left(\sum_{j=1}^{n}j^{d-2}+\sum_{j=n+1}^{\infty}j^{2d-3}\right)
\end{align*}
for some constant $C>0$. Therefore in the case $d<1$, the series for $\mbox{Cov}(Y_n, Z_{n+h})$ converges  as $n\to \infty$. 
By the same argument, the series for $\mbox{Cov}(Z_n, Y_{n+h})$ converges and hence the series for the product
$\mbox{Cov}(Y_n, Z_{n+h})\mbox{Cov}(Z_n, Y_{n+h})$ converges as $n\rightarrow \infty$ if $d<1.$ So in the case $d<1$, 
\begin{align}
\mbox{Cov}(Y_nZ_n, Y_{n+h}Z_{n+h})
=\gamma_y(h)\mbox{Cov}(Z_n, Z_{n+h}) +C(n,h) \label{f}
\end{align}
as  $n\to \infty$ for some constant $C(n,h)$ with uniform bound $0<C<\infty$.
As a particular case of (\ref{f}),
\begin{align}\label{fs}
\mbox{Var}(Y_n Z_n)=\mbox{ Var}(Y_n)\mbox{Var}(Z_n)+C(n)  \;\;\;\text{as $n \to \infty$}
\end{align}
for some  constant $C(n)$ with uniform bound $0<C<\infty$.

iii). Next we prove that the non-stationary process $Z_n$ satisfies:
\begin{align}
&\mbox{Corr}(Z_n, Z_{n+h}) \to 1  \;\;\;\text{as $n \to \infty$} \label{z}
\end{align}
 and
\begin{align}
&\mbox{Var}(Z_n)\to \infty  \;\;\;\text{as $n \to \infty$.} \label{v}
\end{align}
We first show that (\ref{z}) holds under the condition (\ref{v}). In fact,
\begin{align}
\mbox{Cov}(Z_n, Z_{n+h})
&=\mbox{Cov}(X_n+X_{n-1}, X_{n+h}+X_{n+h-1}) \nonumber\\
&=\mbox{Cov}(2\sum_{j=1}^{n-1}Y_j+Y_n, 2\sum_{i=1}^{n+h-1}Y_i+Y_{n+h})\nonumber\\
&=4\sum_{i=1}^{n+h-1}\sum_{j=1}^{n-1}\gamma_{y}(i-j)+2\sum_{j=1}^{n-1}
\gamma_{y}(n+h-j)\nonumber\\
&+2\sum_{i=1}^{n+h-1}\gamma_{y}(n-i)+\gamma_{y}(h).\label{covz}
\end{align}
Let $h=0$ in (\ref{covz}), we have 
\begin{align}
&\mbox{Var}(Z_n)=4\sum_{i=1}^{n-1}\sum_{j=1}^{n-1}\gamma_{y}(i-j)+4\sum_{j=1}^{n-1}\gamma_{y}(n-j)+\gamma_{y}(0) \label{vz}
\end{align}
By replacing the $n$ in (\ref{vz}), 
\begin{align*}
\mbox{Var}(Z_{n+h})
&=4\sum_{i=1}^{n+h-1}\sum_{j=1}^{n+h-1}\gamma_{y}(i-j)+4\sum_{j=1}^{n+h-1}\gamma_{y}(n+h-j)+\gamma_{y}(0)\\
&=4\sum_{i=1}^{n-1}\sum_{j=1}^{n-1}\gamma_{y}(i-j)+4\sum_{j=1}^{n-1}\gamma_{y}(n-j)+\gamma_{y}(0)\\
&+8\sum_{i=1}^{n-1}\sum_{j=n}^{n+h-1}\gamma_{y}(i-j)+4\sum_{i=n}^{n+h-1}\sum_{j=n}^{n+h-1}\gamma_{y}(i-j)+4\sum_{j=0}^{h-1}\gamma_{y}(n+j)\\
&=\mbox{Var}(Z_n)+C(n,h)
\end{align*}
for some constant $C(n,h)$ with uniform bound $0<C<\infty$ since $Y_n$ is a short memory process.
Therefore, 
\begin{align}\label{product}
&\mbox{Var}(Z_n) \mbox{Var}(Z_{n+h})=\mbox{Var}^{2}(Z_n)\left[1+\frac{C(n,h)}{Var(Z_n)}\right].
\end{align}
On the other hand, by  (\ref{covz}) and (\ref{vz}),  
\begin{align}
&\mbox{Cov}(Z_n, Z_{n+h})\nonumber\\
&=4\sum_{i=1}^{n-1}\sum_{j=1}^{n-1}\gamma_{y}(i-j)+4\sum_{j=1}^{n-1}\gamma_{y}(n-j)+\gamma_{y}(0) \nonumber \\
&+4\sum_{i=n}^{n+h-1}\sum_{j=1}^{n-1}\gamma_{y}(i-j)-2\sum_{i=n-h}^{n-1}\gamma_{y}(n-i)\nonumber \\
&+2\sum_{i=0}^{h-1}\gamma_{y}(n+i)+2\sum_{i=n}^{n+h-1}\gamma_{y}(n-i)+\gamma_{y}(h)-\gamma_{y}(0)\nonumber\\
&=\mbox{Var}(Z_n)+4\sum_{i=n}^{n+h-1}\sum_{j=1}^{n-1}\gamma_{y}(i-j)+2\sum_{i=0}^{h-1}\gamma_{y}(n+i)-\gamma_{y}(h)+\gamma_{y}(0) \nonumber\\
&=\mbox{Var}(Z_n)+C(n,h),\label{covz2} 
\end{align}
 where $C(n,h)$ is bounded uniformly by some constant $C>0$. 

Provided that (\ref{v}) holds, i.e., $\mbox{Var}(Z_n)\to \infty$, by (\ref{product}) and (\ref{covz2}),  
\begin{align*}
&\mbox{Corr}(Z_n, Z_{n+h})=\frac{\mbox{Cov}(Z_n, Z_{n+h})}{\sqrt{\mbox{Var}(Z_n)\mbox{Var}(Z_{n+h})}}\\
&=\frac{\mbox{Var}(Z_n)+C(n,h)}{\mbox{Var}(Z_n)\sqrt{[1+\frac{C(n,h)}{\mbox{Var}(Z_n)}]}}\\
&\to 1 \;\;\; \text{as $n \to \infty$},
\end{align*}
which proves that (\ref{z}) is true. 

Now we show that the second property (\ref{v}) of the non-stationary process $Z_n$  also holds.
Since $Y_n \sim$ FARIMA$(0,d-1,0)$ with $1/2< d<1$, we have $\gamma_y(h)<0$ for all $h>0$ and the spectral density at frequency zero
\begin{align}
&f(0)=\gamma_y(0)+2\sum_{h=1}^{\infty} \gamma_y(h)=0. \label{fy0}
\end{align}
For computational convenience, we assume $Y_n \sim$  FARIMA$(0,d,0)$ with $-1/2< d<0$ in the following process. 
Brockwell and Davis (1987) gave the autocovariance function and autocorrelation function for FARIMA$(0,d,0)$ processes with $-1/2<d<1/2$, 
\begin{align}
&\gamma_{y}(0)=\frac{\sigma^2 \Gamma(1-2d)}{\Gamma^{2}(1-d)} \nonumber
\end{align}
and
\begin{align}
&\rho_{y}(h)=\frac{\Gamma(h+d)\Gamma(1-d)}{\Gamma(h-d+1)\Gamma(d)},\;\;\;h\in \mathbb{N},\nonumber
\end{align}
where $\sigma^2=1$ is the variance of innovation $\varepsilon$.
Thus,
\begin{align}
&\gamma_{y}(h)=\frac{\Gamma(1-2d) \sigma^2}{\Gamma(1-d)\Gamma(d)} \frac{\Gamma(h+d)}{\Gamma(h+1-d)}. \label{autocovy}
\end{align}
To prove (\ref{v}), it suffices to prove that the first quantity in (\ref{vz}) goes to infinity as $n\rightarrow \infty$ since $Y_n$ is a short memory process in the covariance sense and therefore the second and third terms in (\ref{vz}) are bounded. 
By collecting terms, we get
\begin{align}
&\sum_{i=1}^{n-1}\sum_{j=1}^{n-1}\gamma_{y}(i-j) =(n-1)\gamma_{y}(0)+2\sum_{h=2}^{n-1}(n-h)\gamma_{y}(h-1)\nonumber \\
&=(n-1)[\gamma_{y}(0)+2\sum_{h=1}^{n-2}\frac{n-1-h}{n-1}\gamma_{y}(h)]\nonumber \\
&=(n-1)[\gamma_{y}(0)+2\sum_{h=1}^{n-2}\gamma_{y}(h)-2\sum_{h=1}^{n-2}\frac{h}{n-1}\gamma_{y}(h)]\nonumber \\
&=(n-1)[-2\sum_{h=n-1}^{\infty}\gamma_{y}(h)-2\sum_{h=1}^{n-2}\frac{h}{n-1}\gamma_{y}(h)]\label{fy1} \\
&=-2(n-1)\sum_{h=n-1}^{\infty}\gamma_{y}(h)-2\sum_{h=1}^{n-2}h \gamma_{y}(h) \label{8t}.
\end{align}
The equality (\ref{fy1}) is from (\ref{fy0}).  Both of the two terms in equation (\ref{8t}) are positive since $\gamma_{y}(h)<0$ if $h>0$. So we prove that (\ref{8t}) goes to infinity by showing that the first term of  (\ref{8t}) goes to infinity as $n\rightarrow \infty$.
\begin{align}
&-2(n-1)\sum_{h=n-1}^{\infty}\gamma_{y}(h)\nonumber \\
&=-2\sigma^2\frac{\Gamma(1-2d) }{\Gamma(1-d)\Gamma(d)}(n-1)\sum_{h=n-1}^{\infty}\frac{\Gamma(h+d)}{\Gamma(h+1-d)}\label{gammah} \\
&\simeq  (n-1)\alpha\sum_{h=n-1}^{\infty}\frac{e^{-h-d}(h+d)^{h+d-1/2}}{e^{-h+d-1}(h+1-d)^{h-d+1/2}}\label{stirling} \\
&=\alpha e^{1-2d}(n-1)\sum_{h=n-1}^{\infty}\left(1+\frac{1-2d}{h+d}\right)^{-(h+d-1/2)}(h+1-d)^{2d-1}\nonumber \\
&\simeq \alpha e^{1-2d}(n-1) e^{2d-1}\sum_{h=n-1}^{\infty}(h+1-d)^{2d-1} \label{limit}\\
&\to \infty \;\;\;\text{as $n \to \infty$ since $2d+1>0.$} \nonumber
\end{align}
In the above equations, $\alpha=-2\sigma^2\frac{\Gamma(1-2d) }{\Gamma(1-d)\Gamma(d)}$;  (\ref{gammah}) is obtained by using (\ref{autocovy}); (\ref{stirling}) is from the Stirling's approximation; and (\ref{limit}) is obtained from the fact that $\lim_{n \to \infty}(1+\frac{1}{n})^{n}=e.$ 
Consequently, (\ref{v}) is true.

From the above analysis, in particular (\ref{f}), (\ref{fs}), (\ref{z}) and (\ref{v}),
\begin{eqnarray}
&&\mbox{Corr}(Y_n Z_n, Y_{n+h} Z_{n+h})\nonumber \\
&=&\frac{\mbox{Cov}(Y_n Z_n, Y_{n+h} Z_{n+h})}{\sqrt{\mbox{Var}(Y_n Z_n)\mbox{Var}(Y_{n+h}Z_{n+h})}}\nonumber \\
&=&\frac{\mbox{Cov}(Y_n, Y_{n+h})\mbox{Cov}( Z_n,Z_{n+h})+C(n,h)}{\sqrt{\mbox{Var}(Y_n)\mbox{Var}(Z_n)+C(n)}\sqrt{\mbox{Var}(Y_{n+h})\mbox{Var}(Z_{n+h})+C(n,h)}} \nonumber \\
&\simeq&\frac{\gamma_y(h)\mbox{Cov}( Z_n,Z_{n+h})}{\sqrt{\gamma_y(0)\mbox{Var}(Z_n)}\sqrt{\gamma_y(0)\mbox{Var}(Z_{n+h})}} \;\;\;\text{since $\mbox{Var}(Z_n), \mbox{Var}(Z_{n+h}) \to \infty$}\nonumber \\
&=& \frac{\gamma_y(h)}{\gamma_y(0)} \frac{\mbox{Cov}(Z_n, Z_{n+h})}{\sqrt{\mbox{Var}(Z_{n})\mbox{Var}(Z_{n+h})}}\nonumber \\
&\to&\mbox{Corr}(Y_n, Y_{n+h}) \nonumber \;\;\;\text{as $n \to \infty$.}
\end{eqnarray}

As $Y_n \sim$ FARIMA$(0,d-1,0)$, $X_n^2 \sim LM(d)$ when $1/2 < d<1$.\\

\noindent \textbf{Acknowledgment}\;\;  The authors thank the referees and the editors for their careful reading of  the manuscript and the valuable comments and suggestions. The authors also thank Micah B. Milinovich for  the suggestion on hypergeometric series in the proof of Theorem \ref{thm3}. 


\begin{thebibliography}{99}  
 
 \bibitem{WB}
Bailey, W. N. (1935).  {\it General Hypergeometric Series.}  Cambridge Tracts in Mathematics and Mathematical Physics. 

\bibitem {BinghamGoldieTeugels}
Bingham, N. H., Goldie, C. M. and Teugels, J. L. (1987).
\textit{Regular Variation}. Cambridge, UK: Cambridge University Press.

\bibitem{PW}
Bondon, P. and  Palma, W. (2007).  A class of antipersistent processes. {\it J. Time Series Anal. } {\bf 28} 261-273.

\bibitem{Brockwell}
Brockwell, P. J. and Davis, R. A.(1987). {\it Time Series: Theory and Methods}. Springer.

\bibitem{Cox}
Cox, D. R. (1977). Contribution to discussion of paper by A. J. Lawrance and N. T. Kottegoda.
{\it J. Roy. Statist. Soc. Ser. A}  {140:34}.

\bibitem{Cramer}
Cram\'{e}r, H. (1946). {\it Mathematical Methods of Statistics}. Princeton University Press, Princeton, NJ.

\bibitem{Dittmann}
Dittmann, I. and  Granger, C. (2002). Properties of nonlinear transformations of fractionally integrated
processes. \textit{J. Econometrics} {\bf 110} 113-133.

\bibitem{FMRT}
Fa\"{y}, G.,  Moulines, E.,  Roueff, F. and  Taqqu, M. S. (2009).   Estimators of long-memory: Fourier versus wavelets. \textit{J. Econometrics} {\bf 151} 159-177. 

\bibitem{Gauss}
Gauss, F. (1866). Disquisitiones generales circa sericm infinitam, {\it Ges. Werke}  {\bf 3}
123-163 and 207-229.

\bibitem{GP}
Geweke, J. and  Porter-Hudak, S. (1983). The estimation and application of long memory
time series models. {\it J. Time Series Anal.}  {\bf 4}  221-238.

\bibitem{GS}
Giraitis, L. and Surgailis, D. (1985). Central limit theorems and other limit theorems for functionals of Gaussian
processes. {\it Probab. Theory  Related Fields} {\bf 70} 191-212.

\bibitem{GJ}
Granger, C. and  Joyeux, R. (1980). An introduction to long memory time series models and fractional differencing.  {\it J. Time Series Anal. } {\bf 1}  15-29.

\bibitem{Guegan}
Gu\a'{e}gan, D. (2005). How can we define the concept of long memory? An econometric survey.  {\it Econometric Rev.}  {\bf 24}  113-149.

\bibitem{HH1}
Ho, H. and  Hsing, T. (1996). On the asymptotic expansion of the empirical processes of long-memory moving averages. {\it Ann. Statist.} \textbf{24}  992-1024.

\bibitem{HH2}
Ho, H. and  Hsing, T. (1997). Limit theorems for functionals of moving averages. {\it Ann. Probab.} \textbf{25}  1636-1669.

\bibitem{HH2}
Ho, H. and  Hsing, T. (2003). A decomposition for generalized U-statistics of long-memory linear processes. {\it Theory and applications of long-range dependence,  Birkh$\ddot{a}$user Boston}, 143-155, Boston, MA.

\bibitem{Honda}
Honda, T. (2000). Nonparametric density estimation for a long-range dependent linear process. {\it Ann. Inst. Statist. Math.} \textbf{52} 599-611. 

\bibitem{Hosking}
Hosking, J. R. M. (1981). Fractional differencing. {\it Biometrika} \textbf{68} 165-176.

\bibitem{Hsing0}
Hsing, T. (1999). On the asymptotic distributions of partial sums of functionals of infinite-variance moving averages. {\it Ann. Probab.}  \textbf{27} 1579-1599.

\bibitem{Hsing}
Hsing, T. (2000). Linear processes, long-range dependence and  asymptotic expansions. {\it Stat. Inference  Stoch. Process.} \textbf{3} 19-29.

\bibitem{HW}
Hsing, T. and Wu, W. B. (2004). On weighted U-statistics for stationary processes. {\it Ann. Statist.} \textbf{32} 1600-1631. 

\bibitem{HDB}
Hurvich, C. M.,  Dro, R. and  Brodsky, J. (1998).  The mean squared error of Geweke and Porter-Hudak's estimator of the memory parameter of a long-memory time-series.  {\it J. Time Series Anal. }  {\bf 19} 19-46.

\bibitem{KT}
Kokoszka, P. S. and  Taqqu, M. S. (1995). Fractional
ARIMA with Stable Innovations. {\it Stochastic Process. Appl.}  \textbf{60} 19-47.

\bibitem{Kulik}
Kulik, R. (2008). Nonparametric deconvolution problem for dependent sequences.  {\it  Electron.  J. Stat.} \textbf{2} 722-740. 
 
\bibitem{Parzen}
Parzen, E. (1981). Time series model identification and prediction variance horizon, In: Findley, ed. {\it Applied Time Series Analysis II}. New York: Academic Press, 415-447.

\bibitem{MH}
Peligrad, M.,  Sang, H., Zhong, Y. and Wu, W. B. (2014). Exact moderate and large deviations for linear processes.  {\it Statist. Sinica} \textbf{24} 957-969. 

\bibitem{MP}
McElroy, T. and Politis, D. (2007). Self-normalization for heavy-tailed time series with long
memory. {\it Statistica Sinica} \textbf{17} 199-220.


\bibitem{Ruppert}
Ruppert, D. (2011). \textit{Statistics and Data Analysis for Financial Engineering}. Springer.

\bibitem{Seneta}
Seneta, E. (1976). \textit{Regularly Varying Functions}. Lecture Notes in Mathematics {\bf 508}, Springer.

\bibitem{ShaoWu}
Shao, X. and Wu, W. B. (2007). Local Whittle estimation of fractional integration for nonlinear processes. {\it Econometric Theory} \textbf{23} 899-929. 

\bibitem{Taqqu}
Taqqu, M. S. (1979). Convergence of integrated processes of arbitrary Hermite rank. {\it Probab. Theory and Related Fields} {\bf 50} 53-83.    

\bibitem{Velasco1}
Velasco, C. (1999 a). Non-stationary log-periodogram regression. {\it J. Econometrics} {\bf 91}
325-371.  

\bibitem{Velasco2}
Velasco, C. (1999 b). Gaussian semiparametric estimation of non-stationary time series. {\it J. Econometrics} {\bf 20} 87-127.  

\bibitem{WuM}
Wu, W. B. and  Mielniczuk, J. (2002). Kernel density estimation for linear processes.  {\it Ann. Statist} \textbf{30} 1441-1459. 


\bibitem{Wu1}
Wu, W. B. (2002). Central limit theorems for functionals of linear processes and their applications. {\it Statist. Sinica}  \textbf{12} 635 - 649.

\bibitem{Wu1'}
Wu, W. B. (2003). Empirical processes of long-memory sequences. {\it Bernoulli}  \textbf{9}(5), 809 - 831.

\bibitem{Wu2}
Wu, W. B. (2006). Unit root testing for functionals of linear processes.
{\it Econometric Theory}  \textbf{22} 1-14.

\bibitem {WuZhao}
Wu, W. B. and  Zhao, Z. (2008). Moderate deviations for
stationary processes.  {\it Statist. Sinica} \textbf{18} 769--782.

\end{thebibliography}
\end{document}